\documentstyle{amsppt}
\magnification 1000 \hsize=6.1in \vsize=8.75in
\input amstex
\TagsOnRight

\document
 

\def\H{{\frak H}}
\def\F{{\frak F}}
\def\C{{\Bbb C}}
\def\R{{\Bbb R}}
\def\Z{{\Bbb Z}}

\def\G{\Gamma}
\def\g{\gamma}

\def\ee{\varepsilon}

\def\ca{{\frak a}}
\def\cb{{\frak b}}

\def\ci{{\infty}}

\def\sa{{\sigma_\frak a}}
\def\sb{{\sigma_\frak b}}

\def\p{\endproclaim \flushpar {\bf Proof: }}

\def\bs{\ \ $\qed$ \vskip 3mm}

\def\Im{{\text{\rm Im}}}
\def\Re{{\text{\rm Re}}}

\topmatter
\title{The classification of higher-order cusp forms}
\endtitle
\author{NIKOLAOS DIAMANTIS and
DAVID SIM (University of Nottingham) } 
\endauthor
\date  nikolaos.diamantis\@nottingham.ac.uk
\enddate
\endtopmatter

\NoRunningHeads
$$\text{\bf 1. Introduction}$$
In this paper we establish the complete classification of cusp forms of 
all orders. Higher-order cusp forms constitute a natural extension of the 
notion of classical cusp forms and they have attracted the 
interest of several researchers in the broader area of automorphic 
forms during the last few 
years. Reasons for this 
interest include its relevance for new approaches to problems related to 
the distribution of modular symbols ([CDO]), to GL$(2)$ $L$-functions 
([DKMO], [F]), to percolation theory ([KZ]) and, more recently,
to Manin's non-commutative modular symbols ([M]) via the connection of 
higher-order forms with iterated integrals ([DS]). This approach has 
already yielded striking successes, e.g. the proof that modular symbols 
have a normal 
distribution ([PR]), the establishment of higher order Kronecker limit 
formulas ([JO]) etc.

The first step towards a classification of spaces of higher-order cusp 
forms was taken in [DO], where the case of order $2$ was settled. 
Unexpectedly, the classification of higher weights proved to be far from 
routine and it 
seems that ultimately this is related to difficulties to identify the 
underlying cohomology in orders higher than $2$. Indeed, we have not found 
yet an Eichler-Shimura-type theorem similar to that proved in [DO]. 
However, in this work we succeed in computing the dimensions of the 
spaces of higher-order forms and in constructing explicit bases that are 
fully computable.

The first two sections deal with technical preliminaries some of which 
have independent interest (e.g. the growth estimates for antiderivatives 
of higher-order cusp forms). The basic theorem, establishing the 
analytic continuation and growth properties of the generalized Poincar\'e 
series on which the basis elements is built, is proved in section 3.2. It 
is based on 
a quite large scale induction 
which relies crucially on the spectral analysis of a certain 
Poincar\'e-type series.
The basis elements are constructed in 
Sections 3.3 and 3.4, where it is also 
proved that they are actually a basis of the space of weight $2$ 
higher-order cusp forms. In section 4 we contruct bases for 
higher-order cusp forms of higher weight. 

$$\text{\bf 2. Definitions and basic estimates}$$
Let $\G \subset PSL_2(\R)$ be a Fuchsian group of the first kind acting 
on the upper half plane $\H$ with non compact quotient $\G\backslash\H$ 
of genus $g$. We assume that there are $m \ge 2$ inequivalent cusps.
As usual we write $x+iy=z\in 
\H$ and $\H^*$ for $\H$ together with the cusps. Let $d\mu z$ be 
the hyperbolic volume form $dx dy/y^2$ and $$V=\int_{\G\backslash\H}
\frac{dx dy}{y^2}$$ 
the volume of $\G\backslash\H$. For a fundamental domain $\F$ fix 
representatives of the inequivalent cusps in $\overline \F$ and give 
them labels such as 
$\ca,\cb$. Use the corresponding scaling matrices $\sa,\sb$ to give 
convenient local coordinates near these cusps as in \cite{I1}, Ch 2. 
The subgroup $\G_\ca$ is the set of elements of $\G$ fixing 
$\ca$ and 
$$
\sa^{-1} \G_\ca \sa= \G_\ci=
\left\{ \pm \left(\smallmatrix 1 & m \\ 0 & 1
\endsmallmatrix\right)
\; \big | \; \ m\in {\Bbb Z}\right\}.
$$
The slash operator $|_k$ defines an action of $PSL_2(\R)$ on 
functions $f:\H \mapsto \C$ by
$$
(f|_k \gamma)(z)=f(\gamma z) (c z+d)^{-k} 
$$
with $\gamma= \left ( \smallmatrix  * & * \\  c & d
 \endsmallmatrix \right ) \in PSL_2(\R)$.
Extend the action to $\Bbb C [PSL_2(\Bbb R)]$ by linearity. We set $j(\g 
,z)= cz+d$.

We now define the space $S_k^t(\Gamma)$ of \it cusp
forms of order $t$ and weight $k \ge 0$ \rm recursively by setting: 
\newline
(i) $S_k^0(\Gamma)=\{0\}$ and \newline
(ii) for $t \ge 1$, by letting $S_k^t(\Gamma)$ be the space of holomorphic 
functions $f: \H \to \C$ such that:

1. $f|_k (\g-1) \in S_k^{t-1}(\G)$, for all $\g \in \G$,

2. $f|_k\pi=f$, for all parabolic $\pi \in \G$ and 

3. for each cusp $\ca$, $(f|_k \sa)(z) \ll e^{-cy}$ as $y \to \ci$ 
uniformly in $x$ for some constant $c>0$ (``vanishing at the cusps").

When the group is clear, we will be using $S_k^t$ instead of $S_k^t(\G)$.

A useful reformulation of this definition, essentially proved in [DKMO] 
is that a holomorphic $f: \H \to \C$ is a $t$-th order cusp form if and only if
it is invariant under the parabolic elements, it satisfies
$$f|_k (\g_1-1) \dots (\g_t-1)=0 \quad \text{for all $\g_i \in \G$}$$ 
and, for each $\g \in \G$ and each cusp $\ca$, $(f|_k \g\sa)(z) \ll 
e^{-cy}$ as $y \to \ci$ 
uniformly in $x$ with $c>0$ and the implied constant depending on $\g$.

Further, by relaxing the third condition to include functions such that, 
for each cusp $\ca$, $(f|_k \sa)(z) \ll y^c$ as $y \to \ci$ 
uniformly in $x$ for some constant $c$, we obtain the space of 
$t$-th order modular forms. We denote it by $M_k^t$.

%
%
%
The next lemma is stated in greater generality than 
what we need it for in 
the sequel because we want it to cover other situations that have arisen 
in our work. 

We define recursively 
certain sets of maps denoted by $\Cal H_t$. First, for 
convenience we use the superscripts $^+$ and $^-$ to indicate absence
and presence of complex conjugation respectively.
We set \newline
(i) $\Cal H_0=\Cal H_1=\{0\}$ and \newline
(ii) for $t >1$, we set 
$$\Cal H_t:=\bigoplus \Sb r+s \le t \\ 1 \le r, s \le t-1 \endSb 
H_r \otimes S_2^s $$
where $H_r \subset Maps(\G, \Bbb C)$ is generated by maps of the form
$$\g \to \prod_{i=1}^m \Big( \int_{z_i}^{\g z_i} f_i(w)dw \Big )^{\pm}$$
for some $z_i \in \H^*$ 
with $f_i \in S_2^{r_i}$ such that $f_i|_2(\cdot-1) \in \Cal H_{r_i}$ and 
$\sum_{i=1}^m r_i=r$.
(Here and in the sequel we will take sums whose upper limit is smaller 
than the lower to be equal to $0.$)
We identify each of these spaces with their images in $Maps(\G, 
S_2^{t-1})$ under the natural projection.

With this notation we can now state a proposition that gives important 
estimates for derivatives and anti-derivatives of weight $2$ cusp forms of 
all orders.
 
\proclaim{Lemma 2.1} For all $t \ge 0$, if $f \in S^t_2$ satisfies 
$f|_2(\cdot -1) \in \Cal H_t$ then, for any cusp $\ca$,
$$\align
& (i) \, \, \Im(z)|f(z)| \ll \sum_{i=0}^{t-1}|\log(\Im(z))|^i
\\
& (ii) \int_{z_0}^z f(w)\, dw \ll \sum_{j=0}^{t}|\log(\Im (\sa^{-1} z))|^j
\quad \text{and}
\\
& (iii) \int_{z_0}^{\g z_0} f(w)\, dw \ll \sum \Sb i+j \le t \\ 0 \le i, j 
\endSb |\log(\Im (\sa^{-1} z))|^i |\log(\Im (\sa^{-1} \g z))|^j
\endalign
$$
for all $z \in \H$. The implied constants are independent of $z$ and 
$\g$.
\p
We will use induction on $t$. For $t=0$, it is trivial. 
Let now $t>0$. Assume the result holds for orders $<t$ and let $f \in 
S^t_2$ satisfy $f|_2(\cdot -1) \in \Cal H_t$. 
To prove (i), let $F_{\ci}$ be the strip of $(x, y)$ with $y>0$ 
and $|x| \le 1/2$ and $F$
the fundamental domain consisting of $z \in F_{\ci}$ such that $|j(\g, z)|>1$
for all $\g \in \G-\G_{\ci}$. Then $\Im(z)|f(z)| \ll 1$ in $F$ 
because $f$ has exponential decay at each cusp. If, on the other hand, 
$z \in F_{\ci}-F$, then
there is $\g \in \G-\G_{\ci}$ and $w \in F$ such that $z=\g w$. 
According to Lemma 1.25 of [Sh], there is a $r>0$ depending only on $\G$ 
such 
that Im$(w) \le 1/(r^2 \Im (z))$. On the other hand,
since $w \in F$, Im$(z)<$Im$(w)$. Therefore, $\log($Im$(z))<\log($Im$(w))
\le -2\log(r)-\log(\Im(z))$
and hence 
$$|\log(\Im(w))| \le 2|\log(r)|+|\log(\Im(z))|. \tag 2.1$$ 
Now, if 
$$(f|_2(\g-1))(z)=\sum \, ^{'} \prod_{i=1}^m(\int_{z_i}^{\g z_i} 
h_i(w) dw)^{\pm} h_0(z) \tag 2.2$$
where the prime indicates that the summation is over $(h_1, \dots, h_m, 
h_0) \in S_2^{l_1} \times \dots \times S_2^{l_m} \times S_2^{l_0}$ 
($l_1+\dots+l_m+l_0 \le t$) with $h_i|_2(\cdot-1) \in \Cal H_{l_i}$ ($i>0$),
then 
$$|\text{Im}(z) f(z)|=|\text{Im}(w) (f|_2 \g)(w)| \le |\text{Im}(w) f(w)|+
\sum \, ^{'}  |\prod_{i=1}^m \int_{z_i}^{\g z_i} h_i(w) dw|
|\text{Im}(w) h_0(w)|.$$
With the boundedness of $\Im(z) |f(z)|$ in $F$ and the inductive 
hypothesis, we deduce that
$$|\text{Im}(z) f(z)| \ll 
1+ \sum \, ^{'} \prod_{n=1}^m \Big ( \sum \Sb i+j \le \l_n \\ 0 \le i, j 
\endSb |\log(\Im ( w))|^i|\log(\Im (\g w))|^j \Big )
\Big (\sum_{i=0}^{l_0-1}|\log(\text{Im}(w))|^i \Big ) 
$$
with the implied constant independent of $w$ and $\gamma$. 
With (2.1) this implies
$$
|\Im(z)f(z)| \ll
1+ \sum \, ^{'} \prod_{n=1}^m \Big ( \sum \Sb i+j \le \l_n \\ 0 \le i, j  
\endSb |\log(\Im (z))|^{i+j} \Big )
\Big (\sum_{i=0}^{l_0-1}|\log(\text{Im}(z))|^i \Big ) \ll
\sum_{i=0}^{t-1} |\log(\Im(z))|^i
$$
for $z \in F_{\ci}-F$ and thus for all $z \in \H$ because 
both sides of (i) are translation invariant. 

To prove (ii), we first note that
$$\int_{z_0}^{\sa z} f(w)\, dw =\int_{\sa^{-1} z_0}^z 
(f|_2 \sa)(w)\, dw.$$
Since $f|_2\sa(\sa^{-1}\g \sa-1)=f|_2(\g-1)\sa$, an inductive argument 
implies that 
$f|_2 \sa \in S_2^t(\G')$ with $\G':=
\sigma_\frak a^{-1}\Gamma \sigma_\frak a$. In fact, a similar argument
 implies that $f|_2 \sa (\cdot -1) \in \Cal H_t(\G')$. Therefore, we can use
(i) to deduce 
$$\Im(z)|(f|_2\sa) (z)| \ll \sum_{i=0}^{t-1}|\log(\Im(z))|^i$$
Further the invariance under the parabolic elements gives
$$
\int_{z}^{z+1} (f|_2 \sa)(w)\, dw=0.
$$
Therefore, 
$$\int_{z_0}^{\sa z} f(w)\, dw \ll 
\int_{\sigma_\frak a^{-1} z_0}^{x+iy}
\frac{\log(\text{Im} w)^{t-1}+\dots+1}{\text{Im}(w)} dw
$$
with $y=$Im$z$ and $x \equiv $Re$z \mod 1$ and $0 \le x <1.$
The last integral equals 
$$\int_{\sigma_\frak a^{-1} z_0}^{x+i\text{Im}(\sa^{-1}z_0)}+
\int_{x+i\text{Im}(\sa^{-1}z_0)}^{x+iy}
\frac{\log(\text{Im} w)^{t-1}+\dots+1}{\text{Im}(w)} dw
$$
where the integration path is, in both cases, a straight segment.
This sum in turn equals:
$$\frac{\log(\text{Im}(\sa^{-1}z_0))^{t-1}+\dots+1}{\text{Im}(\sa^{-1}z_0)} 
(x-\text{Re}(\sa^{-1}z_0))+i\int_{\text{Im}(\sa^{-1}z_0)}^y 
\frac{\log(s)^{t-1}+\dots+1}{s} ds \ll \sum_{j=0}^t 
(\log(\text{Im}z))^j
$$
Replacing $z$ by $\sigma_\frak a^{-1}z$ completes the proof of (ii).

Finally, by employing (2.2) it is easy to see that
$$\int_{z_0}^{\g z_0} f(w) dw=
\int_{z_0}^{\g z} f(w) dw-
\int_{z_0}^{z} f(w) dw+ 
\sum \, ^{'}  \prod_{i=1}^m(\int_{z_i}^{\g z_i}
h_i(w) dw)^{\pm} \int_{z}^{z_0} h_0(w)dw \tag 2.3$$
where the summation is over $(h_1, \dots, h_m,  
h_0) \in S_2^{l_1} \times \dots \times S_2^{l_m} \times S_2^{l_0},$
($l_1+\dots+l_m+l_0 \le t$) with $h_i|_2(\cdot-1) \in \Cal H_{l_i}$,
($i>0$) and the inequality follows from (ii) and the inductive hypothesis.
\bs

$$\text{\bf 3. Bases for $S^t_2(\G)$}$$
$$\text{\bf 3.1 Preliminaries}$$
We collect here some notation and results we will be using frequently in 
the sequel.

Let $C^\infty(\G\backslash\H,k)$ denote the space of 
smooth functions $\psi$ on $\H$ that transform as
$$
\psi(\g z)=\ee(\g,z)^k\psi(z)
$$
for $\g$ in $\G$ and $\ee(\g,z)=j(\g,z)/|j(\g,z)|$. Note that this
notion of weight in general differs from the previous definition of 
weight. 

We define the \it Maass raising and lowering operators \rm by
$$
R_k=2iy\frac{d}{dz} +\frac{k}{2}, \ L_k=-2iy\frac{d}{d\bar z}
-\frac{k}{2}.
$$
It is an elementary exercise to show that
$$
R_k: C^\infty(\G\backslash\H,k)\rightarrow C^\infty(\G\backslash\H,k+2),\ 
L_k:
C^\infty(\G\backslash\H,k)\rightarrow C^\infty(\G\backslash\H,k-2).
$$
For $n>0,$ we write $R^n$ for $R_{k+2n-2} \cdots R_{k+2} R_k$ and 
$L^n$ for $L_{k-2n+2} \cdots L_{k-2}L_k$. 
We also let $L^0$ and $R^0$ be the identity operator.

A very useful fact proved in [JO] (Lemma 9.2) is that if $\g 
\in$PSL$_2(\R)$ and
$$\mu(s, k, F):=F(\gamma) \Im(\g z)^s e(m \g z)\ee(\g, z)^{-k}$$
then
$$\align
R_k\mu(s, k, F)&=2i\mu(s+1, k+2, \frac{d}{dz}F)+(s+\frac{k}{2})\mu(s, k+2, 
F)-4\pi m \mu(s+1, k+2, F) \\
L_k\mu(s, k, F)&=-2i\mu(s+1, k-2, \frac{d}{d 
\bar z}F)+(s-\frac{k}{2})\mu(s, k-2, F)
\tag 3.1 
\endalign
$$
The hyperbolic Laplacian $\Delta=-4y^2 \ d/dz \ d/d\overline{z}$ 
can be realized as $\Delta =-L_{2}R_0  = -R_{-2}L_0.$
Further, for $\tau \in PSL_2(\R)$, we let the operator $\theta_{\tau, k}: 
C^\infty (\G \backslash \H ,k) \rightarrow C^\infty (\tau^{-1}\G \tau 
\backslash
\H ,k)$ be defined by 
$$
\theta_{\tau, k} \psi(z) = \frac{\psi(\tau z)}{\ee(\tau, z)^k}. 
$$
It is easy to verify that this action commutes with the 
raising and lowering operators:
$$
\align
\theta_{\tau, k-2} L_k &= L_k \theta_{\tau, k}, \\
\theta_{\tau, k+2} R_k &= R_k \theta_{\tau, k}. 
\endalign
$$
In stating our bounds, the notation 
$y_\G(z)=\max_\ca (\max_{\g \in \G}(\Im( \sa^{-1}\g z)))$ will often be 
useful.
For example, if $|\psi|$ is smooth with weight $0$ then 
we write
$
\psi(z) \ll y_\G(z)^A,
$
instead of $\psi(\sa z) \ll y^A$ for each cusp $\ca$ as $y 
\rightarrow \infty$. 
We also use the notation
$
y_{\F}(z)=\max_\ca (\Im( \sa^{-1} z))
$
for $z$ in a fundamental domain $\F$. 

We next recall that the usual non-holomorphic Eisenstein series
$$E_{\ca}(z,s)=\sum_{\g \in \G_\ca \backslash \G}   \Im(\sa^{-1} \g z)^s 
$$
is absolutely convergent for $s$ with $\Re(s)>1$ (and uniformly 
convergent for $s$ in compact sets there) and that it has the Fourier 
expansion at the cusp $\cb$
$$
\align
E_{\ca}(\sb z,s)&=\delta_{\ca \cb}y^s +\phi_{\ca \cb}(s)y^{1-s}+\sum_{m 
\neq 0} \phi_{\ca \cb}(m,s)W_s(mz) \\
&=\delta_{\ca \cb}y^s +\phi_{\ca \cb}(s)y^{1-s}+O(e^{-2\pi y}) \tag 3.2
\endalign
$$
as $y \rightarrow \infty$ with an implied constant depending only on $s$ 
and $\G$. 

The hyperbolic Laplacian operates on $L^2(\G \backslash \H)$ the space 
of smooth, automorphic, square integrable functions. Any element 
$\xi$ of $L^2(\G \backslash \H)$ may be 
expanded according to the discrete and continuous spectrum of $\Delta$
(Roelcke-Selberg decomposition):
$$
\xi(z)=\sum_{j=0}^\infty\langle \xi,\eta_j\rangle
\eta_j(z)+\frac{1}{4\pi }\sum_\cb \int_{-\infty}^{\infty}\langle
\xi ,E_\cb(\cdot,1/2+ir)\rangle E_\cb(z,1/2+ir)\,dr,\tag 3.3
$$
where $\{\eta_j\}$ denotes a complete orthonormal basis of Maass
forms, with corresponding eigenvalues $\lambda_j=s_j(1-s_j)$,
which forms the discrete spectrum. As always, we
will write $s_j=\sigma_j+it_j$, chosen so that $\sigma_j \geqslant
1/2$ and $t_j \geqslant 0$, and we enumerate the eigenvalues,
counted with multiplicity, by $0=\lambda_0< \lambda_1 \leqslant
\lambda_2 \leqslant \cdots $. Weyl's law ((11.3) of 
\cite{I1}) implies
$$
\#\{j| \ |\lambda_j| \leqslant T\} \ll T. \tag 3.4
$$

The decomposition (3.3) is absolutely convergent for each fixed $z$ and
uniform on compact subsets of $\H$, provided $\xi$ and $\Delta
\xi$ are smooth and bounded (see, for example, Th. 4.7 and
Th. 7.3 of \cite{I1}). 

For each $j$, the Fourier expansion
of $\eta_j$ is
$$
\eta_j(\sa z)=\rho_{\ca j}(0)y^{1-s_j}+\sum_{m\neq 0}
\rho_{\ca j}(m)W_{s_j}(mz). \tag 3.5
$$
For all but finitely many of the $j$ (corresponding to
$\lambda_{j} < 1/4$) we have $\sigma_j=1/2$ and $\rho_{\ca
j}(0)=0$.   
The constant $\delta_\G$ used throughout 
this paper is chosen so that $1-\delta_\G> \sigma_1 \geqslant 1/2$.

With this notation we now state
\proclaim{Lemma 3.1} For all $z \in \H$, $T \in \Bbb R$ and $n \in \Bbb 
Z_+$
we have 
$$ \align
(i) \, \, & R^n\bigl(\eta_j(z)\bigr), L^n\bigl(\eta_j(z)\bigr) \ll_{\G, n} 
(|t_j|^n +1) y_\G(z)^{1/2}+ (|t_j|^{2n+5} +1) y_\G(z)^{-3/2} 
\\
(ii) \, \, & E_{\ca}(z,1/2+ir) \operatornamewithlimits{\ll}_{\G, T}  
y_\G(z)^{1/2} \quad \text{for all} \, \, r \in [T, T+1]
\\
(iii) \, \, & \int_T^{T+1}  \left| R^n E_{\ca}( z, 1/2+ir)\right|^2 \, dr, 
\ \ \int_T^{T+1}  \left| L^n E_{\ca}( z, 1/2+ir)\right|^2 \, dr \ll 
T^{4n+12} y_\G(z)
\endalign
$$
\p
For a proof of see Lemma 8.2, (11.12) and Lemma 8.3 of [DO] respectively.
\bs

We now define the basic auxiliary functions we will be using in the 
sequel and prove their basic properties.
For $k \in 2\Z$, we consider
$$U_{\ca m}(z,s,k)=\sum_{\g \in \G_\ca \backslash \G}  \Im(\sa^{-1} \g 
z)^s e(m \sa^{-1} \g z) \ee(\sa^{-1} \g, z)^{-k} 
$$
and for simplicity we set $U_{\ca m}(z,s):=U_{\ca m}(z,s,0)$. 
By a direct computation based on (3.1)
$$
\align
R_k U_{\ca m}(z,s,k)&=(s+k/2)U_{\ca m}(z,s,k+2)-4\pi m U_{\ca
m}(z,s+1,k+2)\tag 3.6 \\
L_k U_{\ca m}(z,s,k)&=(s-k/2)U_{\ca m}(z,s,k-2). \tag 3.7
\endalign
$$
\proclaim{Proposition 3.2} 
For $k \in 2\Z$, $U_{\ca m}(z,s,k)$ has a meromorphic 
continuation to all $s$ with Re$(s)>1-\delta_{\G}$ its only pole appearing 
at $s=1$ when $m=k=0$. It is simple with residue $1/V$. Furthermore,
$$U_{\ca 0}(z,s,k) \ll y_{\G}(z)^{\sigma} \quad \text{and} \, \, \,  
U_{\ca m}(z,s,k) \ll y_{\G}(z)^{1/2} \quad (m>0)$$
for these $s$ with the implied constant depending on $s, m, k, \G$.
\p
This is the content of Propositions B and C of [DO]. \bs
Given this analytic continuation we set
$$
P_{\ca m}(z)_2:=y^{-1}U_{\ca m}(z,1,2),
$$
These series are holomorphic for $m>0$ and span $S_2(\G)$ (cf. \cite{JO} 
Th. 3.2). When $m=0$ they satisfy
$$
j(\sb,z)^{-2}P_{\ca 0}(\sb z)_2=\delta_{\ca \cb}-\frac{1}{y V}+O(e^{-2\pi 
y}) \quad \text{as} \, \, y\rightarrow \ci \tag 3.8
$$
and
$$ y^2\frac d{d\overline{z}}P_{\ca 0}(z)_2=\frac{i}{2 V}. \tag 3.9
$$

For $f \in S_2^t$ and $n \in \Bbb Z_{\ge 1}$ we set:
$$
I_{\ca n}(z_n)=\int_{i\ci}^{z_n} \cdots
\int_{i\ci}^{z_2}\int_{i\ci}^{z_1} f_\ca(z_0) dz_0 dz_1 \cdots dz_{n-1}
$$
where $f_\ca(z)=f(\sa z)/j(\sa,z)^2$. For $n \le 0$, we set
$I_{\ca n}(z)=f^{(-n)}_{\ca}(z)$. We observe that $I_{\ca 1}(\sa^{-1}z)=
\int_{\ca}^z f(w)dw$. 
With this notation, we set, for each $r \in \Z_{\ge 0}$ 
$$Q_{\ca m}(z, s, n, r; \bar f)=\sum_{\g \in \G_{\ca} \backslash \G}
\overline{I_{\ca n}(\sa^{-1} \g z)} \, \, \, \Im(\sa^{-1} \g z)^{s} 
e(m \sa^{-1} \g z) \ee(\sa^{-1} \g, z)^{-r}.$$

In this section we will give the domains of initial convergence and 
bounds of this series and its derivative. Their analytic continuation 
will be discussed in the next section.

\proclaim{Proposition 3.3} Let $f \in S_2^t$ be such that $f|_2(\cdot-1) \in
\Cal H_t$. For $k \in 2\Z$ and $\sigma=\Re (s) >1$ the series 
$Q_{\ca m}(z,s, 1, k, \bar f)$ and $Q'_{\ca m}(z,s, 1, k, \bar f)$ 
converge absolutely and uniformly on compacta to analytic functions of $s$. 
For these $s$ and for all $m \ge 0$, $Q_{\ca m}(z,s, 1, k; \bar f)$ 
and $(|m|+1)^{-1} y Q'_{\ca m}(z,s, 1, k; \bar f)$
are bounded by a constant times $y_{\G}(z)^{1/2-\sigma/2}.$
The implied constants are independent of $z$ and $m$.
\p
Let $f \in S_2^t$ be such that $f|_2(\cdot-1) \in
\Cal H_t$. By Lemma 2.1 and the elementary inequality $|\log y| 
<\epsilon (y^{\epsilon}+y^{-\epsilon})$ we have
$$ \int_{\ca}^{\sb z}f(w)dw \ll \sum_{j=-t}^t y^{j\epsilon} \tag 3.10$$ 
for all 
$z$ in $\H$ and hence 
$$ \int_{\ca}^{\g \sb z}f(w)dw=\int_{\ca}^{\sa 
\sa^{-1} \g \sb z} f(w)dw \ll \sum_{j=-t}^t \Im (\sa^{-1} \g \sb 
z)^{j\epsilon}$$ for any cusp $\cb$ and any $z$ in $\H$. The implied 
constant depends solely on $\ca,$ $\epsilon, f$ and $\G$. Further, the 
Fourier expansion of $\int_{\cb}^{\sb z}f(w)dw $ yields 
$$ \int_{\ca}^{\sb z}f(w)dw = \int_\ca^\cb f(w)\, dw +\frac{1}{2\pi i} 
\sum_{n=1}^\infty \frac{a_\cb(n)}{n} e(nz) \tag 3.11 $$ 
with $a_\cb(n)$ the $n$-th Fourier coefficient of $f$ at the cusp $\cb$. 
Thus, 
$$ \int_{\ca}^{\sa z}f(w)dw \ll e^{-2\pi y} \text{ \ as } y \rightarrow 
\infty. \tag 3.12 $$ 
Consequently, since $|e(m \sa^{-1} \g z) \ee(\sa^{-1} \g, z)^{-k}| \le 1$,
$$ \align & Q_{\ca m}(\sa z, s, 1, k; \bar f) \ll \sum_{\g \in 
\G_\ca \backslash \G} |\int_{\ca}^{\g \sa z}f(w)dw | \Im(\sa^{-1} \g \sa 
z)^\sigma \\ & \ll e^{-2\pi y} y^\sigma + \sum \Sb \g \in \G_\ca 
\backslash \G \\ \g \neq I \endSb \left( \sum_{j=-t}^t \Im(\sa^{-1} \g \sa 
z)^{\sigma+j\epsilon}\right) \ll y^{1-\sigma +t\epsilon} \endalign $$ 
for $\sigma> 1+t\epsilon $ as $y \rightarrow \infty$ by (3.2). When $\ca 
\neq \cb$ we have 
$$ Q_{\ca m}(\sb z, s, 1, k; \bar f) \ll \sum_{j=-t}^t E_\ca(\sb z, 
\sigma+j\epsilon) \ll y^{1-\sigma +t\epsilon}$$ 
for $\sigma >1+t\epsilon$ as $y \rightarrow \infty$. Choose $\epsilon 
=(\sigma -1)/2t$ 
for simplicity and we have demonstrated that $$ Q_{\ca m}(z, s, 1, k;
\bar f) \ll y_\G(z)^{1/2-\sigma/2} $$ for $\sigma > 1$ and an 
implied constant depending on $\sigma, \ca, f$ and $\G$ alone. This proves 
the statement for $Q_{\ca m}(z, s, 1, k; \bar f)$.

With (3.1) we deduce that 
$$
\multline
2iy Q'_{\ca m}(z, s, 1, k; \bar f)=
(s+\frac{k}{2}) Q_{\ca m}(z, s, 1, k+2; \bar f)-\\
4\pi m Q_{\ca m}(z, s+1, 1, k+2; \bar f)-
\frac{k}{2} Q_{\ca m}(z, s, 1, k; \bar f)
\endmultline
$$
This together with the continuation and bounds of $Q_{\ca m}$ that we 
have just proved yields the desired result about $Q'_{\ca m}$.
\bs

$$\text{\bf 3.2. The basic theorem}$$

We are ready to state the theorem that will enable us to construct 
the basis elements for $S_2^t$.

We first construct a family of elementary functions in $S_2^t.$ Fix a cusp 
$\ca$
and let $\{f_1, \dots, f_g\}$ be an orthonormal basis of
$S_2$.
For $i_j \in \{1, \dots, g\}$ we set
$$F_{i_1, \dots, i_{t}}(z)=f_{i_1}(z) \int_{\ca}^z f_{i_2}(w)
\Big (\int_{\ca}^w f_{i_3}\dots \Big )dw.$$
It is easy to see with (2.3) and an inductive argument that
$$F_{i_1, \dots, i_{t}}|_2(\g-1)=\sum_{r=1}^{t-1}F_{i_1, \dots, i_r}
\int_{\ca}^{\g \ca} F_{i_{r+1}, \dots, i_t}(w) dw. \tag 3.13$$
It is straightforward, by an inductive argument, to see
that they are invariant under
the parabolic elements, that they  vanish at the cusps and
that $F_{i_1, \dots, i_t}|_2(\g-1) \in S_2^{t-1}$. Hence
$F_{i_1, \dots, i_t} \in S_2^{t}$.

Since these functions will play a fundamental role in the sequel,
we set
$$A_t=\{F_{i_1, \dots, i_t}; i_j \in \{1, \dots, g\}\}.$$
It is clear that, if $f \in A_t$, $f|_2(\cdot-1) \in \Cal H_{t}$, so the
results of the previous sections apply to the elements of this set.

To generate further higher-order cusp forms, we need some functions that
depend on standard cusp forms in a less elementary way than the $F_{i_1,
\dots, i_t}$'s do. Specifically, let $m \ge 0$ and $k \in 2 \Bbb Z$. For 
$f \in S^t_2(\G)$ we set
$$
Z_{\ca m}(z, s, 1, k; \bar f):= 
\sum_{\g \in \G_\ca \backslash \G} 
\overline{\left (\int_{\ca}^{\g \ca} f(w) dw \right )} \Im(\sa^{-1} \g z)^s  
e(m \sa^{-1} \g z) \ee(\sa^{-1}\g, z)^{-k} $$

This function is essentially a generalization of the function $Z_{\ca
m}(z, s; f)$ which was crucial for the construction of a basis 
of the space of second-order cusp forms in [DO]. However, the multiplier 
is slightly modified. The advantage is that in this way we avoid the 
introduction of the function $G_{\ca m}$ used in [DO].

We need to meromorphically continue $Z_{\ca m}$ to a
region that contains $1$. The proof has many similarities to that of the 
corresponding result in [DO].

\proclaim{Theorem 3.4} 
For $f \in A_t$, 
$Z_{\ca m}(z,s, 1, k; \bar f)$ admits a meromorphic continuation to 
$\Re(s)> 1-\delta_\G$. The only possible pole is $s=1$ and 
it can only occur when $k \le 0$. For $k=0$, it is simple.
For $k \ge 2$, and $m \ne 0$, $Z_{\ca m}(z, s, 1, k; \bar f) 
\ll y_{\F}(z)^{1/2}$. 
For $k \ge 2$, $Z_{\ca 0}(z, s, 1, k; \bar f) 
\ll y_{\F}(z)^{\sigma}$. 
For $k=2$, $Z_{\ca 0}(z, 1, 1, 2; \bar f) 
\ll y_{\F}(z)^{1/2}.$ 
The implied constants are independent of $z$ in all cases. 
\p 
We will prove the theorem by induction on $t$. For $t=0$, it is trivial.
Let $t>0$ and suppose, the statement is true for orders $<t$.
If $f \in A_t$, then (2.3) implies
$$\multline
Z_{\ca m}(z,s, 1, k; \bar f)=Q_{\ca m}(z,s, 1, k; \bar f)-
\overline {(\int_{\ca}^{z} f(w) dw)} \, \, U_{\ca m}(z,s, k)+\\
\sum \, ^{'} 
Z_{\ca m}(z, s, 1, k; \bar h) \overline{(\int_{z}^{\ca}h_1(w)dw)} 
\endmultline 
\tag 3.14
$$
where the prime indicates that the summation is over some pairs $(h, h_1)
\in A_r \times A_{t-r}$, $1 \le r \le t-1$.

The meromorphic continuation and bounds of $U_{\ca m}(z, s, k)$ and
$Z_{\ca m}(z, s, 1, k; \bar h)$ in (3.14) are known by Prop. 3.2 and the
inductive hypothesis respectively. Therefore, we only need to
meromorphically continue and bound $Q_{\ca m}(z, s, 1, k; \bar f)$. To
this end we first need to study $Q_{\ca m}(z,s+n+1,-n, k, \bar f)$ for $n
\ge 0$.

\proclaim{Proposition 3.5} Suppose that $f \in A^t$. 
For $k \in \Bbb Z$ and $-n \leqslant 0$ the series $Q_{\ca 
m}(z,s+n+1,-n, k, \bar f)$ has a meromorphic continuation to 
$\Re (s) >1-\delta_\G$. For $k \ge 0$, it is analytic.
Also for these $s$ and $k \ge 0$  we have $
Q_{\ca m}(z,s+n+1,-n, k; \bar f)
\ll e^{-\pi y_\G(z)} $
with the implied constant depending on $n, m, f, k, s$ and $\G$ alone. 
\p
We begin with the formula
$$ 
\overline{g^{(n)}(\g z)}=(-2i)^{-n} \Im(\g z)^{-n-1} \sum_{r=0}^n 
(-1)^{n-r} \ee(\g, z)^{-2r-2}\binom{n}{r} \frac{(n+1)!}{(r+1)!} 
L^r\left(\theta_{\g, -2}(y \overline{g(z)})\right)
$$
valid for any holomorphic $g: \H \to \C$ and for all  $\g$ in $\G$ (see 
\cite{CO} for a proof). 
Set $\G':=\sa^{-1}\G \sa$. Then $\sa^{-1}\G_\ca \sa=\G_\ci$ and
$$
\multline
\ee(\sa, z)^{-k}Q_{\ca m}(\sa z,s,-n, k; \bar f)=
\sum_{\g' \in \G_\ci \backslash \G'}  
\overline{f_\ca^{(n)}( \g' z)} \Im( \g' z)^s e(m  \g' z) \ee(\g', z)^{-k}=\\
(-2i)^{-n} \sum_{r=0}^n (-1)^{n-r} \binom{n}{r} \frac{(n+1)!}{(r+1)!} 
\sum_{\g' \in \G_\ci \backslash \G'} 
L^r(\theta_{\g', -2}(y \overline{f_\ca(z)})) 
\Im(\g' z)^{s-n-1} e(m  \g' z) \ee(\g', z)^{-(2r+2)-k}
\endmultline
$$
Now, if $f|_2(\g-1)=\sum \, ^{'}(\int_{\ca}^{\g \ca} h(w) dw ) h_1$ 
($\g \in \G$), with $(h, h_1) \in A_{r}(\G) \times A_{t-r}(\G)$, ($r \ge 
1$), then 
$f_{\ca}|_2(\g'-1)=\sum \, ^{'}(\int_{\sa ^{-1} \ca}^{\g' \sa^{-1} \ca} 
h_{\ca}(w) dw ) (h_1)_{\ca}$ for $\g' \in \G'$ and
$(h_{\ca}, (h_1)_{\ca}) \in A_r(\G') \times A_{t-r}(\G').$ Therefore,
$$
L^r\left(\theta_{\g', -2}(y \overline{f_\ca(z)})\right)= L^r \left(y 
\overline{(f_\ca|_2\g')(z)}\right)=
L^r \left(y 
\overline{f_\ca(z)}\right)+ \sum \, ^{'}
\overline{\int_{\sa^{-1} \ca}^{\g' \sa^{-1} \ca} h_{\ca}(w) dw } L^r \left(y 
\overline{(h_1)_\ca(z)}\right) 
$$ 
It is also easy to see (cf. Prop. E of [DO]) that 
for all $G: \H \to \C$, $$L^r \left(y \overline{G_\ca(z)}\right)=L^r \left .\left(y 
\overline{G(z)}\right)\right|_{\sa z} \ee(\sa, z)^{2r+2}$$
so
$$
\multline
\sum_{\g' \in \G_\ci \backslash \G'}
L^r\left(\theta_{\g', -2}(y \overline{f_\ca(z)})\right) 
\Im(\g' z)^{s-n-1} 
e(m  \g' z) \ee(\g', z)^{-(2r+2)-k}=\\
L^r(y \overline{f(z)})|_{\sa z} 
U_{\ca m}( \sa z, s-n-1, 2r+2+k) \ee(\sa, z)^{-k} +\\
\sum \, ^{'}\Big [L^r(y \overline{h_1(z)}) \times \\
\Big (\sum_{\g \in \G_{\ca} \backslash \G} 
\Big (\overline{\int_{\ca}^{\g \ca} h(w)dw} \Big )
\Im(\sa^{-1} \g z)^{s-n-1} e(m \sa^{-1} \g z) \ee(\sa^{-1} \g, 
z)^{-(2r+2)-k} \Big ) \Big ]|_{\sa z} \ee(\sa, z)^{-k}\\
=L^r(y \overline{f(z)})|_{\sa z} 
U_{\ca m}( \sa z, s-n-1, 2r+2+k) \ee(\sa, z)^{-k}+\\
\sum \, ^{'} L^r(y \overline{h_1(z)})|_{\sa z}
Z_{\ca m} (\sa z, s-n-1, 1, 2r+2+k; \bar h)\ee(\sa, 
z)^{-k}
\endmultline
\tag 3.15$$
Since $h \in A_{r}$, $r <t$, the inductive hypothesis implies that 
$Z_{\ca m}$ has a 
meromorphic continuation for Re$(s-n-1)>1-\delta_{\G}$. By Prop. 3.2 the same 
holds for $U_{\ca m}$. Moreover, since $2r+2 \ge 2$, by the inductive hypothesis
and Prop. 3.2 we deduce that we obtain an analytic function when $k \ge 0$.

The function $L^r\left(y \overline{f(z)}\right)$ (and $L^r\left(y 
\overline{h(z)}\right)$) has exponential decay at 
every cusp $\cb$ because 
$$ 
 \theta_{\sb, -2r-2} L^r \left(y \overline{f(z)}\right)  
= L^r \left(\theta_{\sb, -2} y \overline{f(z)}\right)
= L^r \left(y \overline{ j(\sb z)^{-2} f(\sb z)}\right)
= L^r \left(y \sum_{n=1}^\infty \overline{ a_\cb (n)e(nz)}\right).
$$
Hence 
$$
L^r \left(y \overline{f(z)}\right) \ll y_\G(z)^{r+1} e^{-2\pi y_\G(z)} 
\tag 3.16
$$
for an implied constant depending on $r, f$ and $\G$. Therefore, with 
(3.15), (3.16), Prop. 3.2
and the inductive 
hypothesis, we have for $\Re(s)>1-\delta_\G$ and $k \ge 0$:
$$
Q_{\ca m}(z,s+n+1,-n, k; \bar f)  \ll  e^{-\pi y_\G(z)}. 
$$
\bs


We are now ready to prove the analytic continuation of $Q_{\ca m}(z,s, 1,
k, \bar f)$.

\proclaim{Proposition 3.6} Let $m \ge 0$ and $k \in 2 \Bbb Z$. For $f 
\in A_t$ the series 
$Q_{\ca m}(z,s, 1, k; \bar f)$ 
has continuation to a meromorphic function of $s$ with $\Re (s) 
>1-\delta_\G$. 
For $k>0$, we obtain an analytic function. For $k=0$, it has only 
a simple pole at $s=1$ with residue 
$\frac{-\overline{a_{\ca}(m)} }{2\pi i m }$ if $m \ne 0$ and $0$ 
otherwise.
For $k \ne 0$ or $k=m=0$, we have 
$Q_{\ca m}(z,s, 1, k; \bar f) \ll  y_\G(z)^{1/2}$. 
For $k=0$, $m \ne 0$, 
$(s-1)Q_{\ca m}(z,s, 1, 0; \bar f) \ll  y_\G(z)^{1/2}.$
The implied constants depend on $s$, $m$, $f$ and $\G$. 
\p 
We first prove the result for $k=0$. By Proposition 3.3, $Q_{\ca 
m}(z,s, 1, 0; \bar f)$ is square integrable 
for $\Re(s)>1$ and the spectral decomposition yields
$$
\multline
Q_{\ca m}(z,s, 1, 0; \bar f)=
\sum_{j=0}^\infty\langle Q_{\ca 
m}(\cdot,s, 1, 0; \bar f),\eta_j\rangle
\eta_j\\
+\frac{1}{4\pi }\sum_\cb \int_{-\infty}^{\infty}\langle
Q_{\ca m}(\cdot,s, 1, 0; \bar f) ,E_\cb(\cdot,1/2+ir)\rangle 
E_\cb(z,1/2+ir)\,dr.
\endmultline \tag 3.17
$$
We recall the Prop. 9.3 and Cor. 9.4 of [JO]:
\proclaim{Lemma 3.7} 
Let $\xi_1$, $\xi_2$ 
and $\psi$ be any smooth $\G$ invariant functions 
(not necessarily in $L^2(\G \backslash \H)$). If $(\Delta - 
\lambda) \xi_1 =\xi_2$, $(\Delta - \lambda') \psi =0$ and
$$
\align
\xi_1, R_0\xi_1, \Delta \xi_1 & \ll y_\G(z)^A,
\\ \psi, R_0 \psi & \ll y_\G(z)^B
\endalign
$$
for $A+B<0$ and $R_0=2iy \frac{d}{dz}$ the raising operator, then
$$
\langle \xi_1,\psi \rangle = \frac 1{\lambda' -\lambda} \langle 
\xi_2,\psi \rangle.
$$
\endproclaim
We will apply this lemma to
$\xi_1 =Q_{\ca m}(z,s,n, 0; \bar f)$ ($n \in \Bbb Z$) and 
$\psi = \eta_j$. (3.1) implies that for all $n \in \Z$,
$$
\align
(\Delta - s(1-s))Q_{\ca m}(z,s, n, 0; \bar f)  = & -8 \pi i m Q_{\ca 
m}(z,s+2,n-1, 0; \bar f)\\
&+ 4\pi ms Q_{\ca m}(z,s+1,n, 0; \bar f)
+ 2i s Q_{\ca m}(z,s+1,n-1, 0; \bar f).
\endalign
$$
Next, we have $\eta_j(z)$, $R_0\eta_j(z) \ll y_\G(z)^{1/2}$ by Lemma 
3.1(i) and 
$$
Q_{\ca m}(z,s,1, 0; \bar f), R_0Q_{\ca m}(z,s, 1, 0; \bar f), 
\Delta Q_{\ca m}(z,s, 1, 0; \bar f) \ll y_\G(z)^{1/2-\sigma/2}
$$
for $\sigma= \Re(s)>1$ by Proposition 3.3 and Proposition 3.5. 
So we may use Lemma 3.7 to get, for $\Re(s)>2$,
$$
\multline
\langle Q_{\ca m}(\cdot,s,1, 0; \bar f),\eta_j\rangle = 
\frac{1}{(s_j-s)(1-s_j-s)}
\Big( 
-8 \pi i m \langle Q_{\ca m}(\cdot,s+2,0, 0; \bar f),\eta_j\rangle
\\
+ 4\pi ms \langle Q_{\ca m}(\cdot,s+1,1, 0; \bar f),\eta_j\rangle
+ 2i s \langle Q_{\ca m}(\cdot,s+1, 0, 0; \bar f),\eta_j\rangle \Big).
\endmultline
$$
We can repeat this procedure $W$ times in all to obtain, again for 
$\Re(s)>2$,
$$
\langle Q_{\ca m}(\cdot,s,1, 0; \bar f),\eta_j\rangle = \sum_l 
\frac{P_l(m,s)}{R_l(s_j,s)} 
\langle Q_{\ca m}(\cdot,s+W+c_l,1-d_l, 0; \bar f),\eta_j\rangle, \tag 3.18
$$
with integers $c_l,d_l$ satisfying $0 \leqslant c_l,d_l \leqslant W$, 
$d_l \le W+c_l$, $P_l(m,s)$ a polynomial in $m,s$ alone of degree $W$ in 
$m$ and of 
degree $W$ in $s$ and $R_l(s_j,s)$ a polynomial 
in $s_j,s$ alone of degree $2W$ in $s_j$ and of degree $2W$ in $s$. In fact
$$
R_l(s_j,s)=\prod_b (s_j-b-s)(1-s_j-b-s) \tag 3.19
$$
where, for each $l$, the 
product is over some subset of integers $b$ in 
$\{0,1, \cdots ,2W\}$ of cardinality $W$.

If $d_l=0$, then we have
$$
Q_{\ca m}(z, s+W+c_l,1, 0 ; \bar f) \ll y_\G(z)^{1/4-W/2}\tag 3.20
$$
by Proposition 3.3, for $W \geqslant 1$. Hence
$$
\langle Q_{\ca m}(\cdot,s+W+c_l, 1, 0; \bar f),\eta_j\rangle 
\ll \sqrt{ || y_\G(z)^{-1/4}|| \cdot || \eta_j ||} = 
\sqrt{ || y_\G(z)^{-1/4}|| } \ll 1. \tag 3.21
$$
For $0< d_l \leqslant W$, Prop. 3.5 implies that
$$Q_{\ca m}(z,s+W+c_l,1-d_l, 0; \bar f) \ll  e^{-\pi y_\G(z)}$$
and hence $\langle Q_{\ca m}(z,s+W+c_l,1-d_l, 0; \bar f),\eta_j\rangle 
\ll  1.$
Therefore, for $j>0$,
$
\langle Q_{\ca m}(\cdot,s,1, 0; \bar f),\eta_j\rangle
$
is an analytic function of $s$ for $\Re(s)>1 - \delta_\G$ and satisfies
$$
\langle Q_{\ca m}(\cdot,s,1, 0; \bar f),\eta_j\rangle \ll |s_j|^{-2W} 
\ll |\lambda_j|^{-W} \tag 3.22
$$
for implied constants depending on $s,m,W,f$ and $\G$ alone and with the 
dependence on $s$ being uniform on compacta.

For $j>0$ and for all $n \ge 0$ we can now use (3.4), Lemma 
3.1 (i) and 
(3.22) to get 
$$\sum_{T\leqslant |\lambda_j|< T+1} <Q_{\ca m}(\cdot, s, 1, 0; 
\bar f), 
\eta_j> R^n\bigl(\eta_j(z) \bigr) \ll T^{1-W}((T^{n/2}+1) y_\G(z)^{1/2}+ 
(T^{n+7/2} +1) y_\G(z)^{-3/2}).$$
Hence for $W=6+n$,
$$\sum_{j=1}^{\infty} <Q_{\ca m}(\cdot, s, 1, 0; \bar f), \eta_j> 
R^n (\eta_j(z)) \ll y_\G (z)^{1/2} \tag 3.23
$$
for all $s$ with Re$(s)>1-\delta_{\G}$. 
The sum converges uniformly for $s$ in compact sets 
with $\Re (s) > 1-\delta_\G$ giving an analytic function of $s$. 
When $n>0$, $R^n$ eliminates $\langle Q_{\ca
m}(\cdot,s, 1, 0; \bar f),\eta_0\rangle \eta_0(z)$.

For $j=0$, the constant eigenfunction is $\eta_0=V^{-1/2}$. If $m \ne 0$, 
by unfolding we obtain
$$
\langle Q_{\ca m}(\cdot,s, 1, 0;\overline{f}),
\eta_0\rangle \eta_0 =  \frac{-\overline{a_{\ca}(m)} \ 
\Gamma(s-1)}{2\pi i m (4\pi m)^{s-1}} 
= \frac{-\overline{a_{\ca}(m)} }{2\pi i m }\left( \frac{1}{s-1} +O(1)\right)
$$
as $s\rightarrow 1$ 
for $f_\ca(z)=
\sum_{m=1}^\infty a_{\ca}(m) e(mz).$
If $m=0$, the same process gives $0$.

With arguments similar to those used for 
the discrete spectrum we now consider the continuous spectrum. For $P_l, 
R_l, c_l$ and $d_l$ identical to $(3.18)$, Lemma 3.1 (ii) gives
$$ \multline
\langle Q_{\ca m}(\cdot,s,1, 0; \bar f),E_\cb(\cdot, 1/2+ir)\rangle = \\
\sum_l \frac{P_l(m,s)}{R_l(1/2+ir,s)} 
\langle Q_{\ca m}(\cdot,s+W+c_l,1-d_l, 0; \bar f),E_\cb(\cdot, 1/2+ir)
\rangle, 
\endmultline
\tag 3.24
$$
which is true for $\Re(s)> 2$ initially. 

With (3.19), (3.20), Prop. 3.5 and Lemma 3.1 (ii) we see 
that (for $W\geqslant 1$) the right side of $(3.24)$ 
converges and gives the analytic 
continuation of the left side to $\Re(s)>1-\delta_\G$.
 Now, for all $n \ge 0$,
$$
\multline
\int_T^{T+1} \langle Q_{\ca m}(\cdot,s,1, 0; \bar f),
E_\cb(\cdot, 1/2+ir)\rangle R^n E_\cb(z_0, 1/2+ir)\, dr\\
 =\sum_l P_l(m,s) \int_T^{T+1}  \frac{\langle 
Q_{\ca m}(\cdot,s+W+c_l,1-d_l, 0; \bar f),
E_\cb(\cdot, 1/2+ir)\rangle}{R_l(1/2+ir,s)}  
R^n E_\cb(z_0, 1/2+ir)\, dr \\
 =\sum_l P_l(m,s) \int_T^{T+1} \int_{\F} \frac{Q_{\ca m}(z,s+W+c_l,1-d_l, 
0; \bar f)}{R_l(1/2+ir,s)} \overline{E_\cb(z, 1/2+ir)} R^n E_\cb(z_0, 1/2+ir)\, 
d\mu z \, dr. \endmultline\tag 3.25 
$$ 
The integrand satisfies 
$$ 
\multline \frac{Q_{\ca m}(z,s+W+c_l,1-d_l, 0; \bar f)}{R_l(1/2+ir,s)} 
\overline{E_\cb(z, 1/2+ir)} R^n E_\cb(z_0, 1/2+ir) \ll \\ 
|r|^{-2W+n}y_\G(z)^{1/4-W/2}y_\G(z)^{1/2}y_\G(z_0)^{1/2} \endmultline $$ 
by (3.19), (3.20), Prop. 3.5, Lemma 3.1(ii) and the easily proved 
identity $$R^n E_\ca(z,s) =s(s+1) \cdots (s+n-1) U_{\ca 0}(z,s,2n).$$ 
Thus 
the double integral in (3.25) is absolutely and uniformly convergent and 
we may interchange the limits of integration to obtain $$ \sum_l P_l(m,s)  
\int_{\F} Q_{\ca m}(z,s+W+c_l,1-d_l, 0; \bar f) 
\int_T^{T+1}\frac{\overline{E_\cb(z, 1/2+ir)}}{R_l(1/2+ir,s)} R^n 
E_\cb(z_0, 1/2+ir)\, dr \, d\mu z. 
\tag 3.26 $$ 
Also 
$$ \multline 
\int_T^{T+1}\frac{\overline{E_\cb(z, 1/2+ir)}}{R_l(1/2+ir,s)} R^n 
E_\cb(z_0, 1/2+ir)\, dr 
\\ \ll T^{-2W} \sqrt{ \int_T^{T+1} |E_\cb( z, 
1/2+ir)|^2\, dr \cdot \int_T^{T+1} |R^n E_\cb( z_0, 1/2+ir)|^2\, dr}. 
\endmultline $$ 
So, with Lemma 3.1 (iii), (3.26) is bounded by a constant 
times 
$$\sum_l |P_l(m,s)| T^{12-2W+2n} \int_{\F} y_\G(z)^{3/4-W/2} \, d\mu z 
\,\cdot y_\G(z_0)^{1/2}. 
$$ 
This 
means that, for $W$ chosen large enough, 
$$ \int_{-\infty}^{\infty}\langle 
Q_{\ca m}(\cdot,s, 1, 0; \bar f), E_\cb(\cdot,1/2+ir)\rangle R^n 
E_\cb(z,1/2+ir)\,dr \ll y_\G(z)^{1/2}. $$ 

To combine the information we have collected for the discrete and the 
continuous part of (3.17) we observe that with (3.23) and preceeding 
discussion we can interchange summation and differentation to get
$$
R^n\bigl(\sum_{j=1}^{\infty} <Q_{\ca m}(\cdot, s, 1, 0; \bar f), \eta_j> 
\eta_j(z) \bigr)=\sum_{j=1}^{\infty} <Q_{\ca 
m}(\cdot, s, 1, 0; \bar f), \eta_j> R^n\bigl(\eta_j(z) \bigr)
$$ 
and that Lemma 10.2 of [DO] implies 
$$ \multline 
\frac{1}{4\pi}\sum_\cb 
R^n \int_{T}^{T+1} \langle Q_{\ca m}(\cdot,s, 1, 0; \bar f) 
,E_\cb(\cdot,1/2+ir)\rangle E_\cb(z,1/2+ir)\, dr= \\ 
\frac{1}{4\pi}\sum_\cb \int_{T}^{T+1} \langle Q_{\ca m}(\cdot,s, 1, 0; 
\bar f) ,E_\cb(\cdot,1/2+ir)\rangle R^n E_\cb(z,1/2+ir)\, dr 
\endmultline $$ 
for all $T \in \R.$ Therefore, 
$$\multline
R^n Q_{\ca m}(z,s, 1, 0; \bar f) =\langle Q_{\ca m}(\cdot,s, 
1, 0; \bar f), \eta_0\rangle R^n\eta_0 +
\sum_{j=1}^\infty\langle Q_{\ca m}(\cdot,s, 1, 0; \bar f),\eta_j\rangle
R^n \bigr (\eta_j(z) \bigl ) \\
+\frac{1}{4\pi }\sum_\cb \int_{-\infty}^{\infty}\langle
Q_{\ca m}(\cdot,s, 1, 0; \bar f) ,E_\cb(\cdot,1/2+ir)\rangle 
R^n E_\cb(z,1/2+ir)\,dr
\endmultline 
$$
and it has a continuation to a meromorphic function of $s$ with 
Re$(s)>1-\delta_{\G}$. For these values, 
$(s-1) Q_{\ca m}(z, s, 1, 0; \bar f)$,  
$R^n Q_{\ca m}(z, s, 1, 0; \bar f)$ ($n>0$) and 
$Q_{\ca 0}(z, s, 1, 0; \bar f)$ are all $\ll y_{\G}(z)^{1/2}$ 
The only pole is at $s=1$ and it comes from 
the contribution of 
$\eta_0$ when $n=0$.

To pass to general $k$'s we apply the operators $R_r$ successively, using 
(3.1), to obtain, for $k>0$:
$$            
\multline
Q_{\ca m}(z,s, 1, k; \bar f)=\frac{1}{s(s+1) \cdots (s+\frac{k}{2}-1)} 
\Big ( R^{k} Q_{\ca m}(z,s, 1, 0; \bar f) \\
+ \sum_{i=-1}^{k-1} \sum_{j=1+i}^k \, ^{'} p_{i, j}(m,s)Q_{\ca m}(z,s+j, 
-i, k; \bar f) \Big ) \endmultline 
\tag 3.27
$$
with polynomials $p_{i, j}$ in $m$ and $s$. Here the prime indicates that 
we exclude the term corresponding to $(i, j)=(-1, 0).$ 
Thanks to Propositions 3.3, 3.5 (for $-i \le 0$) and the meromorphic 
continuation and growth of $R^{k} Q_{\ca m}(z, s, 1, 0; f)$ we just 
proved, the identity (3.27) implies Prop. 3.6. For $k<0$, we work in a 
similar way. \bs

\it End of proof of Th. 3.4 \rm By Propositions 3.2 and 3.6, the inductive 
hypothesis and (3.14), we 
deduce that for $k>0$ $Z_{\ca m}(z, s, 1, k; \bar f)$ is holomorphic in 
$s$ and that for $k=0$ the only possible pole is at $s=1$ which is 
simple.
This completes the proof of the analytic 
continuation of $Z_{\ca m}$.
 
To prove the bounds, we recall from (3.11) and $(3.12)$ that, for $z \in 
\F$, $ \int_{\ca}^{\sa z}f(w)dw \ll 1$ and 
$ \int_{\ca}^{\sb z}f(w)dw \ll  e^{-2\pi y}$, if $\ca \ne \cb$ 
as $y \to \infty$. The desired bound 
follows from this,
Propositions 3.2 and 3.6, the inductive hypothesis and (3.14).
For $m=0$, $k=2$ and $s=1$ we deduce the bound from these inequalities and 
$U_{\ca 0}(\sb z, 1, 2) \ll \delta_{\ca \cb} y+1$, as $y \to \infty$ 
(see (3.8)).
\bs 

$$\text{\bf 3.3 A family of functions of $S_2^t$.}$$

In this section we construct a family of $t$-order cusp forms based on the
analytic continuation of $Z_{\ca m}(z, s, 1, 2; f)$ established in
Section 3.2. The construction is carried out in three steps.

In the first step, since we are mainly interested in the weight
according to $j(\g, z)$ rather than $\ee(\g, z)$, we set
$$\align Z_{\ca m}(z, s; \bar f)& :=y^{-1}Z_{\ca m}(z, s+1, 1, 2; \bar 
f)\\
&=\sum_{\g \in \G_\ca \backslash \G} \overline{\Big (\int_{\ca}^{\g \ca}
f(w) dw \Big )}
\Im(\sa^{-1} \g z)^{s} e(m \sa^{-1} \g z) j(\sa^{-1} \g, z)^{-2}.
\tag 3.28
\endalign
$$
According to Theorem 3.4, $Z_{\ca m}(z, s, \bar F_{i_1, \dots, 
i_{t-1}})$ is analytic for Re$(s)>-\delta_{\G}$.
It is easy to see that
$$\multline
Z_{\ca m}(\cdot, 0, \bar F_{i_1, \dots, i_{t-1}})|_2(\g-1)=\\
\Big (\overline{\int_{\ca}^{\g^{-1} \ca} 
F_{i_1, \dots, i_{t-1}}}\Big ) P_{\ca m}+
\sum_{r=1}^{t-2} \Big (\overline{\int_{\ca}^{\g^{-1} \ca} 
F_{i_1, \dots, i_r}}\Big )Z_{\ca m}(\cdot, 0, 
\bar F_{i_{r+1}, \dots, i_{t-1}}).
\endmultline$$
Further, for Re$(s)$ large we have:
$$
\frac{d}{d\bar z}Z_{\ca m}(z, s; \bar F_{i_1, \dots, i_{t-1}})=
\frac{is}{2y^2}Z_{\ca m}(z, s+1, 1, 0; \bar F_{i_1, \dots, i_{t-1}})
\tag 3.29$$

In the second step, if $f_{i_t}=\sum_l a_l P_{\ca m_l}$, we set, for $j 
\ge 1$,
$$Z_{i_j, \dots, i_t}=\sum_l a_l
Z_{\ca m_l}(\cdot, 0; \bar F_{i_j, \dots, i_{t-1}}).$$
An inductive argument implies that, for $t \ge 2$,
$$Z_{i_1, \dots, i_{t}}|_2(\g-1)=\sum_{r=1}^{t-1}
\Big (\overline{\int_{\ca}^{\g^{-1} \ca} F_{i_1, \dots, i_r}}
\Big ) Z_{i_{r+1}, \dots, i_t} \tag 3.30$$
In view of (3.29), we apply the same linear combination to
Res$_{s=1}Z_{\ca m}(\cdot, s, 1, 0; \bar F_{i_{1}, \dots, i_{t-1}})$ and to  
$\frac{-\overline{a_{\ca}(m)} }{2\pi i m }$, for $m>0$. We denote them by 
$R_{i_1, \dots, i_t}$ and $a_{i_1, \dots, i_t}$ respectively. Then, 
(3.14) and Prop. 3.2 give
$$R_{i_1, \dots, i_t}=a_{i_1, \dots, i_t}+ \sum_{r=1}^{t-2}
\overline{\Big ( \int_{z}^{\ca}F_{i_1, \dots, i_r}\Big)}
R_{i_{r+1}, \dots, i_t}. \tag 3.31$$
For convenience we have set $Z_j:=f_j.$
(3.30) and induction imply that
$$Z_{i_1, \dots, i_{t}}|_2(\g_1-1) \dots (\g_{t-1}-1)=
\overline{\int_{\ca}^{\g^{-1}_{1} \ca} f_{i_1}} \dots
\overline{\int_{\ca}^{\g^{-1}_{t-1} \ca} f_{i_{t-1}}}
f_{i_t}  \tag 3.32$$

In the third step, we suitably modify the functions constructed so far to
obtain holomorphic forms. To state a lemma we will need, we 
recursively define the following functions:
$$S_{i_1}:=\overline{\int_z^{\ca} f_{i_1} (w)dw} \quad \text{and} $$
$$S_{i_1, \dots, i_t}:=\sum_{r=1}^{t}\overline{\int_z^{\ca} F_{i_1, 
\dots, i_r}(w) dw} \, \, S_{i_{r+1}, \dots, i_t}.$$
We also set $S_{i_{k+1}, \dots, i_k}=1$, $S_{i_j, \dots, i_k}=0$ for 
$j>k+1$ and $S_{i_0}=1$. 
\proclaim{Lemma 3.8} Let $t \in \Z_{\ge 2}$. For every 
$i_1, \dots, i_t \in \{1, \dots, g\}$ we have \newline
(i) For $m \ne 0$,
$$R_{i_1, \dots, i_t}=\sum_{j=0}^{t-2} S_{i_1, \dots, i_j} 
\cdot a_{i_{j+1}, \dots, i_t}.$$
(ii) 
$$\text{Res}_{s=1}Z_{\ca 0}(\cdot, s, 1, 0; \bar F_{i_{1}, \dots, 
i_t})=\frac{1}{V}S_{i_1, \dots, i_{t}}.$$
\p
(i) is proved by induction in $t$. By (3.14), it is clear for 
$t=2$.
If the result holds for orders $<t$, then (3.31) implies that
$$R_{i_1, \dots, i_t}=a_{i_1, \dots, i_t}+ \sum_{r=1}^{t-2}
\overline{\Big ( \int_{z}^{\ca}F_{i_1, \dots, i_r}\Big)}
\sum_{j=r}^{t-2}S_{i_{r+1}, \dots, i_j}a_{i_{j+1}, \dots, i_t}.$$
By the definition of $S_{i_j \dots i_k}$ for $j>k$, the 
inner sum can be written
in the form $\sum_{j=1}^{t-2}S_{i_{r+1}, \dots, i_j}a_{i_{j+1}, \dots, 
i_t}$ and thus
$$R_{i_1, \dots, i_t}=a_{i_1, \dots, i_t}+ \sum_{j=1}^{t-2}
\sum_{r=1}^{t-2} 
\Big ( \overline{ \int_{z}^{\ca}F_{i_1, \dots, i_r}}
S_{i_{r+1}, \dots, i_j} \Big )a_{i_{j+1}, \dots, i_t}.$$
Since $S_{i_j, \dots, i_k}=0$ for $j>k+1$, the inner sum equals $S_{i_{1}, 
\dots, i_{j}}.$
This proves the identity for $t$. \newline
(ii) follows from (3.14), Prop. 3.2 and a straightforward induction 
argument because,
as shown in the proof of Prop. 3.6., $a_{i_1, \dots, i_t}$ in (3.14) is 
$0$. \bs
By (3.29) and Lemma 3.8(i), there is a linear combination
$Z(i_1, \dots, i_{t-1})$
of $Z_{i_1, i_2}, \dots, Z_{i_1, \dots, i_{t-1}}$ for $m \ne 0$ such that
$$\frac{d}{d\bar z}(Z_{i_1, \dots, i_t}-Z(i_1, \dots, i_{t-1}))=
\frac{i}{2y^2}(S_{i_1, \dots, i_{t-2}}a_{i_{t-1}, i_t}) \tag 3.33$$
Since
$$\text{Res}_{s=1}(Q_{\ca m}(z, s, 1, 0; \bar F_{i_{t-1}}))=
\frac{-\overline{a_{\ca}(m)}}{2\pi i m}=
2i\overline{<f_{i_{t-1}}, P_{\ca m}(\cdot)>},$$ 
(where $<\cdot , \cdot>$ is the usual Petersson scalar product),
by the definition of $a_{i_{t-1}, i_t}$ we have,
$$a_{i_{t-1}, i_t}=2i\overline{<f_{i_{t-1}}, f_{i_t}>}.$$
Therefore, if $i_{t-1} \ne i_t$, (3.33) and the orthonormality of the 
basis imply
that $Z_{i_1, \dots, i_t}-Z(i_1, \dots, i_{t-1})$ is holomorphic. 
If $i_{t-1}=i_t$, (3.33) again implies that
$$
\Big ( Z_{i_1, \dots, i_t}-Z(i_1, \dots, i_{t-1}) \Big )-
\Big ( Z_{i_1, \dots, i_{t-2}, 1, 1}-Z(i_1, \dots, i_{t-2}, 1) \Big )
$$
is holomorphic.

In order to keep track of the basis elements we will construct in the 
sequel, we indicate the conjugation of $\overline{\int_{i 
\ci}^{\g_j i \ci} f_{i_j}(w)dw}$ by a minus sign in the 
notation of the corresponding subscript. Specifically, in view of (3.32)
we set, for every $i_1, \dots, i_{t} \in \{1, \dots, g\}$ with
$(i_{t-1}, i_t) \ne (1, 1)$,
$$\Cal Z_{-i_1, \dots, -i_{t-1}, i_t}=\cases  (-1)^{t-1} 
\left (Z_{i_1, \dots, i_t}-
Z(i_1, \dots, i_{t-1}) \right ) & \, \text{if} \, \, i_{t-1} \ne i_t \\
 (-1)^{t-1} (Z_{i_1, \dots, i_t}-Z(i_1, \dots, i_{t-1})-
Z_{i_1, \dots, 1, 1}-Z(i_1, \dots, i_{t-2}, 1))
& \, \text{if} \, \, i_{t-1}= i_t.  \endcases $$
The reason we have added the factor $(-1)^{t-1}$ is so that $\Cal Z$ 
satisfy (3.32) without the $\g_j$'s being inverted on the right-hand side.
For $t=1$, we set, for $i<0,$
$$\Cal Z_{i}=\Cal Z_{-i}=f_{i}:=f_{-i}.$$
\proclaim{Theorem 3.9} For $t \in \Z_{\ge 2},$
$i_1, \dots, i_{t-1} \in \{-1, \dots, -g\}$, $i_t \in \{1, \dots, g\}$ 
and $(i_{t-1}, i_t) \ne (-1, 1)$, we have $\Cal Z_{i_1, 
\dots, i_t} \in S_2^t$.
\p
They are holomorphic by construction. The invariance under the 
parabolic elements of all 
$Z_{i_1, \dots, i_r}$'s is deduced by (3.30) and the fact that $F_{i_1, 
\dots, i_r} \in S_2^{t}$. (3.32) implies that 
$$\Cal Z_{i_1, \dots, i_t}|_2(\g_1-1) \dots (\g_{t}-1)=0.$$
The growth condition proved in Theorem 3.4 
implies that each $Z_{i_1, \dots, i_r}$ and thus 
$\Cal Z_{i_1, \dots, i_t}$ is $\ll y_{\F}(z)^{-1/2}$. (Recall 
that $Z_{\ca m}(z, s+1, 1, 2; \bar f)$ is divided by $y$ in 
$Z_{\ca m}(z, s; \bar f)$). Considering 
the Fourier expansion of $\Cal Z_{i_1, \dots, i_t}$ we deduce its 
vanishing at the cusps.

Now, (3.30) and Th. 3.4 imply that $Z_{i_1, \dots, i_r}|_2 \g \ll
 y_{\F}(z)^{-1/2}$ for each $\g \in \G$. Therefore, the same estimate 
holds for the holomorphic $\Cal Z_{i_1, \dots, i_t}|_2 \g$ and a look at 
its Fourier expansion implies the vanishing at the cusps for each $\g$.
 
By (the second formulation of) the definition of $S_2^t$ we deduce the 
result.
\bs

In the next section we will also need a family of functions lying outside 
$S_2^t$. They are the analogue of the function $Z_{\mu, \mu}+2iV 
\overline{<\mu, \mu>} P_{\ca 0}(z)_2$ in Prop. 5.2 of [DO]. 
If $\ca_i$
($i=1, \dots, m$) is a set of inequivalent cusps for $\G$, set
$$f_{g+i}:=P_{\ca_i 0}-P_{\ca_m 0} \qquad \text{for} \, \, i=1, \dots,
m-1.$$
Then $\{f_{g+1}, \dots, f_{g+m-1}\}$ a basis of the space $\Cal E_2$ of 
Eisenstein series of $M_2$ (cf. [GO]). 
For $i_j \in \{1, \dots, g\}$, we set
$$\Cal Z'_{-i_1, \dots, -i_{t-1}, i_t}=\cases  \Cal Z_{-i_1, \dots, i_t}
& \quad \text{if} \, \, i_{t-1} \ne i_t \\
 (-1)^{t-1} \Big ( Z_{i_1, \dots, i_t}- Z(i_1, \dots, i_{t-1})-
2iV \cdot Z_{\ca_m 0}(\cdot, 0 ; \bar F_{i_1, \dots, i_{t-2}})
\Big )
& \quad \text{if} \, \, i_{t-1}= i_t.  \endcases $$
The last term is understood to be $P_{\ca_m 0}(\cdot)_2$ when $t=2$. 
By Lemma 3.8(ii), (3.29) and (3.33) we observe that
$\Cal Z'_{-i_1, \dots, -i_{t-1}, i_t}$ is holomorphic. 
By the holomorphicity, Theorem 3.4 and (3.8), $\Cal Z'_{-i_1, 
\dots, i_t}$ vanishes at all cusps except $\ca_m$ in the case $t=2$.
Again the reason for the factor $(-1)^{t-1}$ is that $\Cal Z'$ now
satisfies (3.32) without the $\g_j$'s being inverted on the right-hand 
side.

$$\text{\bf 3.4 Construction of a basis of $S_2^t$}$$
We will use the functions defined in Section 3.3 to build recursively 
bases for all $S_2^t$'s. We first introduce some notation and prove an 
elementary lemma. 
First, following [R], for an increasing finite sequence $j_1, \dots 
j_{t-1}$ we call a \it shuffle \rm of type $(r, t)$ a pair 
$(\phi, \psi)$ of order-preserving maps
$$
\phi: \{j_1, \dots, j_{r-1}\} \to \{1, \dots, t-1\} \quad \text{and} \, \,
\psi: \{j_{r}, \dots, j_{t-1}\} \to \{1, \dots, t-1\}
$$
whose images are disjoint and complementary. Since the specific underlying 
sequence will be understood in each case, we denote their set 
simply by $\Cal{S}_{r, t}.$
We then have
\proclaim{Lemma 3.10} For $F \in S_2^{r}$, $G \in M_0^{t-r+1}$ we have
$$(F \cdot G)|_2(\g_1-1) \dots (\g_{t-1}-1)=
\sum_{(\phi, \psi) \in \Cal{S}_{r, t}} F|_2(\g_{\phi(1)}-1) \dots 
(\g_{\phi(r-1)}-1) \cdot G|_0(\g_{\psi(r)}-1) \dots (\g_{\psi(t-1)}-1).$$
\p 
By Th. 2.2 of [CD], we have
$$(F \cdot G)|_2(\g-1)=F|_2(\g-1) \cdot G+ F \cdot G|_0(\g-1)
+F|_2(\g-1) \cdot G|_0(\g-1).$$
This means that each time we apply a $\g-1$ ($\g \in \G$) on $F \cdot G$,
either $F$ or $G$ or both are acted upon by $\g-1$ too. Now, $F$ 
(resp. $G$) is annihilated by any products with $r$ (resp. $t-r+1$) 
factors of the form $\g-1$. Therefore, the only non-vanishing terms left 
after 
$(\g_1-1) \dots (\g_{t-1}-1)$ is applied on $F \cdot G$ are the 
products of the form 
$$F|_2(\g_{i_1}-1) \dots (\g_{i_{r-1}}-1) \cdot G|_0(\g_{j_1}-1)\dots 
(\g_{j_{t-r}}-1)$$
with $i_1< \dots < i_{r-1}$, $j_1 < \dots < j_{t-r}$
and 
$i_1,  \dots, i_{r-1}, j_1, \dots, j_{t-r}$ covering $\{1, \dots, 
t-1\}$. In particular, by the last two facts the sets of $i_k$ 
and 
$j_k$ are disjoint. This implies Lemma 3.10.
\bs
Next, for $j \in \{\pm 1, \dots, \pm g \} \cup \{g+1, \dots, g+m-1\}$
we set $<f_j, \g>$ for 
$\int_{i}^{\g i} f_j(w)dw$, when $j>0$ and $\overline{\int_{i 
\ci}^{\g i \ci} f_{-j}(w)dw}$, for $j<0$. 
We also set $\Cal A$ for the space generated by maps $\phi: \G^{l} \to 
S_2$ defined by
$$\phi(\g_1, \dots, \g_l)=<f_{i_1}, \g_1> \dots <f_{i_l}, 
\g_l>f_{i_{l+1}}$$ with $-i_j=i_{j+1}=1$ for at least one $j \in \{1, 
\dots, l\}$. (In section 4 we will need the analogue of this space with 
$f_{i_{l+1}}$ of higher weight. In an effort to simplify notation, $\Cal 
A$ will stand for the space in the case of weight $2$ and when we need it 
for higher weights, we will indicate it by a subscript.) 

We denote by $I'$ the set of vectors $(i_1, \dots, i_t)$ with entries 
in $\{\pm 1, \dots, \pm g\}$ for which there is no $j \in \{1, \dots, 
t-1\}$ such that $-i_j=i_{j+1}=1$. We also let $I$ be the set of 
$(i_1, \dots, i_t) \in I'$, with $i_t>0$.

Starting with $\Cal Z_i$, suppose now that, for each $s<t$, the forms 
$\Cal Z_{i_1, \dots, i_s},$ $(i_1, \dots, i_{s}) \in I$ satisfy
$$\Cal Z_{i_1, \dots, i_s}|_2(\g_1-1) \dots (\g_{s-1}-1)=
<f_{i_1}, \g_1>\dots <f_{i_{s-1}}, \g_{s-1}> f_{i_s}
+\phi(\g_1, \dots, \g_{s-1}) 
\tag 
3.34$$
for some $\phi \in \Cal A$.

We claim that for every 
$(i_1, \dots, i_t) \in I$, 
there are $\Cal Z_{i_1, \dots, i_{t}}$ satisfying (3.34) with $s=t$.

First, with Lemma 3.10, 
$$\multline
[\Cal Z_{i_1} \int_{i}^z \Cal Z_{i_2, \dots, i_t}(w)dw]|_2(\g_1-1) 
\dots (\g_{t-1}-1)=\\
f_{i_1} <f_{i_2}, \g_1> \dots <f_{i_{t-1}}, \g_{t-2}> 
\int_{i}^{\g_{t-1} i} f_{i_t}(w)dw +\phi(\g_1, \dots, \g_{t-1}).
\endmultline
\tag 3.35
$$
Hence, for 
$(i_1, \dots, i_{t-1}) \in I$ and $i_t \in \{1, \dots, g\}$ 
we set 
$$\Cal Z_{i_1, \dots, i_t}:=\Cal Z_{i_t} \int_{i}^z \Cal Z_{i_1, 
\dots, i_{t-1}}(w)dw.$$

Next, if $i_1<0$, $(i_1, i_2) \in I$, and $(i_3, \dots, i_t) \in I$,
$$\multline
[\Cal Z_{i_1, i_2} \int_{i}^z \Cal Z_{i_3, \dots, 
i_t}(w)dw]|_2(\g_1-1) \dots (\g_{t-1}-1)=\\
\sum_{(\phi, \psi) \in \Cal{S}_{2, t}} 
\overline{\int_{i}^{\g_{\phi(1)} i}f_{i_1}(w)dw} 
f_{i_2} \prod_{j=3}^{t-1} <f_{i_j}, \gamma_{\psi(j)}> 
\int_{i}^{\g_{\psi(t)} i} f_{i_{t}}(w)dw+ \phi(\g_1, \dots, 
\g_{t-1})
\endmultline
$$
for some $\phi \in \Cal A$.
By the definition of shuffles, 
$\int_{i}^{\g_{t-1} i} f_{i_t}(w)dw$
appears as the last factor and thus unconjugated in each summand of the 
right-hand side except for that corresponding to $\phi(1)=t-1$. Therefore, 
by 
(3.35), there is a linear combination of 
$\Cal Z_{j_1} \int_{i}^z \Cal Z_{j_2, \dots, j_t}(w)dw$'s 
denoted by $A(f_{i_1}, \dots, f_{i_t})$, such that 
$$\multline
[\Cal Z_{i_1, i_2} \int_{i}^z \Cal Z_{i_3, \dots, 
i_t}(w)dw-A(f_{i_1}, \dots, f_{i_t})]|_2(\g_1-1) \dots (\g_{t-1}-1)=\\
\overline{\int_{i}^{\g_{t-1} i} f_{i_1}(w)dw} 
f_{i_2} <f_{i_3}, \gamma_1> \dots <f_{i_{t-1}}, \g_{t-3}> 
\int_{i}^{\g_{t-2} i} f_{i_t}(w)dw+ \phi(\g_1, \dots, \g_{t-1})
\endmultline
\tag 3.36
$$
for some $\phi \in \Cal A$.
Hence for $(i_{t-1}, i_t), (i_{1}, \dots, i_{t-2}) \in I$ and $i_{t-1}<0$ 
we can set 
$$\Cal Z_{i_1, \dots, i_t}:=\Cal Z_{i_{t-1}, i_t} \int_{i}^z \Cal 
Z_{i_1, \dots, i_{t-2}}(w)dw-A(f_{i_{t-1}}, f_{i_t}, f_{i_1}, \dots, 
f_{i_{t-2}}).$$

Further, if $i_1, i_2 <0$ and $(i_1, i_2, i_3), (i_4, \dots, i_t) \in I$,
$$\multline
[\Cal Z_{i_1, i_2, i_3} \int_{i}^z \Cal Z_{i_4, \dots, 
i_t}(w)dw]|_2(\g_1-1) \dots (\g_{t-1}-1)=\\
\sum_{(\phi, \psi) \in \Cal{S}_{3, t}} 
\overline{\int_{i}^{\g_{\phi(1)} i} f_{i_1}(w)dw}
\overline{\int_{i}^{\g_{\phi(2)} i} f_{i_2}(w)dw}
f_{i_3}\prod_{j=4}^{t-1} <f_{i_j}, \gamma_{\psi(j)}> 
\int_{i}^{\g_{\psi(t)} i} f_{i_{t}}(w)dw+ \\
\phi(\g_1, \dots, \g_{t-1})\endmultline
$$
for some $\phi \in \Cal A$.
The only summand that does not appear on the 
right-hand side of (3.35) and (3.36) for an appropriate permutation of 
$i_1, \dots, i_{t},$ is
$$ \overline{\int_{i}^{\g_{t-2} i} f_{i_1}(w)dw}
\overline{\int_{i}^{\g_{t-1} i} f_{i_2}(w)dw}
f_{i_3}<f_{i_4}, \gamma_1> \dots \int_{i}^{\g_{t-3} i} 
f_{i_t}(w)dw. 
$$
Therefore, by (3.35) and (3.36), there is a linear combination of 
terms $\Cal Z_{j_1} \int_{i}^z \Cal Z_{j_2, \dots, 
j_t}(w)dw$'s and 
$\Cal Z_{j_1, j_2} \int_{i}^z \Cal Z_{j_3, \dots, j_t}(w)dw$'s 
denoted by $B(f_{i_1}, \dots, f_{i_t})$, such that 
$$[\Cal Z_{i_1, i_2, i_3} \int_{i}^z \Cal Z_{i_4, \dots, 
i_t}(w)dw-B(f_{i_1}, \dots, f_{i_t})]|_2(\g_1-1) \dots (\g_{t-1}-1)=$$
$$\overline{\int_{i}^{\g_{t-2} i} f_{i_1}(w)dw}
\overline{\int_{i}^{\g_{t-1} i} f_{i_2}(w)dw}
f_{i_3}<f_{i_4}, \gamma_1> \dots \int_{i}^{\g_{t-3} i} 
f_{i_t}(w)dw+ \phi(\g_1, \dots, \g_{t-1})$$
for some $\phi \in \Cal A$.
Hence for $(i_{t-2}, i_{t-1}, i_t), (i_{1}, \dots, i_{t-3}) \in I$ and 
$i_{t-2}, i_{t-1}<0$ we can set
$$\Cal Z_{i_1, \dots, i_t}:=\Cal Z_{i_{t-2}, i_{t-1}, i_t} \int_{i}^z 
\Cal Z_{i_1, \dots, i_{t-3}}(w)dw-B(f_{i_{t-2}}, f_{i_{t-1}}, f_{i_t}, 
f_{i_1}, \dots, f_{i_{t-3}}).$$

Continuing this way we construct, for all $(i_1, \dots, i_t) \in I$, 
functions $\Cal Z_{i_1, \dots i_t}$ satisfying (3.34).
The last element, whose indices $i_1, 
\dots, i_{t-1}$ are negative, is obtained directly from Th. 3.9.

We also define recursively a similar family of functions involving $\Cal 
Z'$: Starting with $\Cal Z'_i=\Cal Z_i$, suppose that, for each $s<t$,  
$\Cal Z'_{i_1, \dots, i_s},$ $i_j \in \{\pm 1, \dots  \pm g \}$ ($i_s>0$) 
satisfy 
$$\Cal Z'_{i_1, \dots, i_s}|_2(\g_1-1) \dots (\g_{s-1}-1)=
<f_{i_1}, \g_1>\dots <f_{i_{s-1}}, \g_{s-1}> f_{i_s}. \tag 3.34'$$
Then 
$$
[\Cal Z'_{i_1} \int_{i}^z \Cal Z'_{i_2, \dots, i_t}(w)dw]|_2(\g_1-1) 
\dots (\g_{t-1}-1)=
f_{i_1} <f_{i_2}, \g_1> \dots <f_{i_{t-1}}, \g_{t-2}> 
\int_{i}^{\g_{t-1} i} f_{i_t}(w)dw 
$$
and for $i_j \in \{\pm 1, \dots  \pm g \}$ ($i_{t-1}>0$) we set
$$\Cal Z'_{i_1, \dots, i_t}:=\Cal Z'_{i_t} \int_i^z \Cal Z'_{i_1, 
\dots, i_{t-1}}(w)dw.$$
Next, if $i_1<0$, 
$$\multline
[\Cal Z'_{i_1, i_2} \int_{i}^z \Cal Z'_{i_3, \dots, 
i_t}(w)dw]|_2(\g_1-1) \dots (\g_{t-1}-1)=\\
\sum_{(\phi, \psi) \in \Cal{S}_{2, t}} 
\overline{\int_{i}^{\g_{\phi(1)} i}f_{i_1}(w)dw} 
f_{i_2} \prod_{j=3}^{t-1} <f_{i_j}, \gamma_{\psi(j)}> 
\int_{i}^{\g_{\psi(t)} i} f_{i_{t}}(w)dw.
\endmultline
$$
As before, there is a linear combination of 
$\Cal Z'_{j_1} \int_i^z \Cal Z'_{j_2, \dots, j_t}(w)dw$'s 
denoted by $A(f_{i_1}, \dots, f_{i_t})$, such that 
$$\multline
[\Cal Z'_{i_1, i_2} \int_{i}^z \Cal Z'_{i_3, \dots, 
i_t}(w)dw-A(f_{i_1}, \dots, f_{i_t})]|_2(\g_1-1) \dots (\g_{t-1}-1)=\\
\overline{\int_{i}^{\g_{t-1} i} f_{i_1}(w)dw} 
f_{i_2} <f_{i_3}, \gamma_1> \dots <f_{i_{t-1}}, \g_{t-3}> 
\int_{i}^{\g_{t-2} i} f_{i_t}(w)dw
\endmultline
$$
so, if $i_{t-1}<0$ we set 
$$\Cal Z'_{i_1, \dots, i_t}:=\Cal Z'_{i_{t-1}, i_t} \int_{i}^z \Cal 
Z'_{i_1, \dots, i_{t-2}}(w)dw-A(f_{i_{t-1}}, f_{i_t}, f_{i_1}, \dots, 
f_{i_{t-2}}).$$

Continuing in this way we construct, for all $i_j \in \{\pm 1, \dots, \pm 
g\}$ ($i_t>0$), 
functions $\Cal Z'_{i_1, \dots i_t}$ satisfying (3.34), but, this time, 
without a $\phi$.
The last element, whose indices $i_1, 
\dots, i_{t-1}$ are negative, is obtained directly by the construction at 
the end of section 3.3.
  
By construction these functions are holomorphic and satisfy the stated
functional equation. They have at most polynomial growth at each cusp 
$\ca,$ and vanish
at $\ca \neq \ca_m.$ Because of this, they are invariant under 
$\pi_{\ca}$ for $\ca \ne \ca_m.$

To prove that the $\Cal Z$'s span $S_2^t$, we will need a lemma that 
generalizes Prop. 5.2. of [DO]. 
\proclaim{Lemma 3.11} Let $t \ge 3$. Suppose that for some 
$F \in S_2^t$ and $c_{i_3, \dots, i_t} \in \Bbb C$, 
$$F|_2(\g_1-1) \dots (\g_{t-1}-1)=\sum c_{i_3, \dots, 
i_t}\overline{\int_{i}^{\g_1 i} f_1(w)dw}
\int_{i}^{\g_2 i} f_1(w)dw <f_{i_3}, \g_3> \dots f_{i_t}$$
for all $\g_i \in \G$, where the sum ranges over all $i_j \in \{ 
\pm1, \dots, \pm g\}$.
Then all $c_{i_3, \dots, i_t}$ vanish.
\p As shown above, 
$$\Cal Z'_{-1, 1, \dots, i_t}|_2(\g_1-1) \dots (\g_{t-1}-1)=
\overline{\int_{i}^{\g_1 i} f_1(w)dw}
\int_{i}^{\g_2 i} f_1(w)dw <f_{i_3}, \g_3> \dots f_{i_t}
$$
and hence $F-\sum c_{i_3, \dots, i_t} \Cal Z'_{-1, 1, \dots, i_t}$ is 
annihilated by 
$(\g_1-1) \dots (\g_{t-1}-1)$.  
Since $F$ and $\Cal Z'_{-1, 1, \dots}$ are of at most polynomial growth at 
the cusps, for every $\g_i \in \G$, we have
$$(F-\sum c_{i_3, \dots, i_t} \Cal Z'_{-1, 1, \dots, i_t})|_2(\g_1-1) 
\dots 
(\g_{t-2}-1)=\sum_{i=1}^{g+m-1} \chi_i(\g_1, \dots, 
\g_{t-2}) f_i$$ 
for $\chi_i: \G^{t-2} \to \Bbb C$. 
The left-hand side is annihilated upon the application of one more 
$\g-1$ and hence the identity $\g 
\delta-1=(\g-1)(\delta-1)+(\g-1)+(\delta-1)$ 
implies that each $\chi_i$ is a group homomorphism in terms of each 
$\g_j$. 

By Eichler-Shimura isomorphism, 
$$\chi_i(\g_1, \dots, \g_{t-2})=\sum_j \left ( a^i_j(\g_1, \dots, 
\g_{t-3}) 
\int_{i}^{\g_{t-2} i} h_j(w)dw+b^i_j(\g_1, \dots, \g_{t-3}) 
\overline{\int_{i}^{\g_{t-2} i} g_j(w)dw} \right )$$ for some $a^i_j, 
b^i_j: 
\G^{t-3} \to \Bbb C$, $g_j \in S_2$ and $h_j \in M_2$. The injectivity of 
the Eichler-Shimura isomorphism implies that each $a^i_j$, $b^i_j$ is a 
homomorphism on each of the arguments and hence, by induction, we deduce 
that
$$
\multline
(F-\sum c_{i_3, \dots, i_t} \Cal Z'_{-1, 1, \dots, i_t})|(\g_1-1)\dots 
(\g_{t-2}-1)= \\
\sum_{i_1, \dots, i_{t-1}} \lambda_{i_1, \dots, i_{t-1}} 
<f_{i_1}, \g_1> \dots <f_{i_{t-2}}, \g_{t-2}> f_{i_{t-1}}
\endmultline
\tag 3.37
$$ 
with $i_j \in \{\pm 1, \dots, \pm g\} \cup \{g+1, \dots, g+m-1\}$ and
$i_{t-1}>0$. Since $\Cal Z'_{i_1, \dots, i_{t}}$ is invariant under 
$\pi_{\ca_i}$ when $i \ne m$, the identity 
$(\g-1)(\pi-1)=(\g \pi \g^{-1}-1)\g-{\pi-1}$ implies that
$$\Cal Z'_{i_1, \dots, i_{t}}|_2(\g_1-1) \dots (\g_{t-2}-1)$$
vanishes when one of the $\g_j$'s equals $\pi_{\ca_i},$ $i \ne m$.
Since $F \in S_2^t$, this then also holds for both sides of (3.37). 
Therefore, none of the $f_{i_j},$ $j<t-1$ appearing in (3.37) can be 
non-cuspidal. Indeed, for $0<k<m,$  $<f_{i_j}, \pi_{\ca_k}> \ne 0$ iff 
$i_j=g+k$ (recall the definition of $f_{g+1}$). Hence, for all $\g_1, 
\dots, \g_{j-1}, \g_{j+1}, \dots, \g_{t-2} \in 
\G$ and for $\g_j=\pi_{\ca_k},$ ($k=1, \dots, m-1$), the RHS of (3.37) 
equals
$$\multline \sum_{i_1, \dots, i_{j-1}, i_{j+1}, \dots, i_{t-1}}
\lambda_{i_1, \dots, i_{j-1}, g+k, i_{j+1}, \dots i_{t-1}} \times \\
<f_{i_1}, \g_1> \dots <f_{i_{j-1}}, \g_{j-1}>
<f_{g+k}, \pi_{\ca_k}> <f_{i_{j+1}}, \g_{j+1}> \dots <f_{i_{t-2}}, 
\g_{t-2}>f_{i_{t-1}}
=0.
\endmultline$$
Therefore, by the injectivity of Eichler-Shimura 
isomorphism, $\lambda_{i_1, \dots, i_{j-1}, g+k, i_{j+1}, i_{t-1}}=0.$
Hence, only cusp forms $f_j$ appear on the RHS of (3.37), so if $\{\g_1, 
\dots, \g_{t-2}\}$ contains $\pi_{\ca_m}$, both sides of 
(3.37) vanish. Since $F \in S_2^t$, this 
implies that $\sum c_{i_3, \dots, i_t} \Cal Z'_{-1, 1, \dots, 
i_t}|_2(\g_1-1) \dots (\g_{t-2}-1)=0$
if at least one of the $\g_i$'s is $\pi_{\ca_m}.$ However, we can 
show, by induction that, when $i_1<0$, 
$$\Cal Z'_{i_1, i_2, \dots, i_t}|_2(\pi_{\ca_m}-1) \dots (\g_{t-2}-1)=
f_{i_t}<f_{i_{t-1}}, \g_{t-2}>\dots \int_i^{\pi_{\ca_m} i} \Cal Z'_{-i, 
i}(w)dw$$ if $(i_1, i_2)=(-i, i)$ ($i>0$) and $0$ otherwise.
For $t=3$, it is straightforward. If $i_r>0$, $i_{r+1}, \dots, 
i_{t-1}<0$, then 
$\Cal Z'_{-i, i, i_3, \dots}$ is by definition equal to 
$\Cal Z'_{i_{r+1} \dots i_t} \int_i^z \Cal Z'_{-i, i, i_3, \dots, i_r}$ 
minus a linear combination of 
products of the form $\Cal Z'_{i_{r+k}, \dots, i_t} \int_i^z \Cal 
Z'_{j_1, j_2, j_3, \dots}$ ($k>1$) which we denote by $C(f_{-i}, f_i, 
f_{i_3}, \dots)$. The indices $j_1, j_2, \dots$  
in the expression of $C(f_{i_1}, f_{i_2}, \dots)$ 
are obtained by interspersing, in 
their original order, subsets of $\{i_{r+1}, \dots, i_{t-1}\}$ among 
the indices $i_1, i_2, i_3, \dots, i_t$. 
Now, as in Lemma 3.10, we observe that, if one of 
the $\g_i$ is parabolic and $i_{r+k}, \dots i_{t-1}<0$, $i_{r+k-1}>0$ 
then, 
for $k \ge 1$,
$$
\align
[\Cal Z'_{i_{r+k} \dots i_t} &\int_i^z \Cal Z'_{j_1, j_2, \dots, 
i_{r+k-1}}]|_2(\g_1-1)\dots (\g_{t-2}-1)=\\
&\sum_{(\phi, \psi) \in \Cal S}
\Cal Z'_{i_{r+k}, \dots i_t}|_2(\g_{\phi(1)}-1) \dots 
(\g_{\phi(t-r-k)}-1) \times \\
&\int_z^{\g_{\psi(r+k-2)}z}[\Cal Z'_{j_1, j_2, \dots, 
i_{r+k-1}}|_2(\g_{\psi(1)}-1)\dots (\g_{\psi(r+k-3)}-1)](w)dw.
\tag 3.38
\endalign
$$
Therefore, by the second part of the inductive hypothesis, the application 
of $(\pi_{\ca_m}-1) \dots(\g_{t-2}-1)$ on $\Cal Z'_{-i, i, \dots}$ will 
eliminate all terms in $C(f_{-i}, f_{i}, f_{i_3}, \dots)$ except for those 
with $j_1=-i$, $j_2=i$, thus 
obtaining, by (the first part of) the inductive hypothesis, products 
of the form
$$\Big ( \int_i^{\g_{\psi(1)} i}\Cal Z'_{-i, i}(w)dw \Big )  \dots 
<f_{i_{r+k-1}}, \g_{\psi(r+k-2)}> <f_{i_{r+k}}, \g_{\phi(1)}> \dots 
<f_{i_{t-1}}, 
\g_{\phi(t-r-k)}> f_{i_t}.$$
These, upon varying $k$, yield 
$$
\multline
\sum_{(\phi, \psi) \in \Cal S_{r, t-1}}
 \Big ( \int_i^{\g_{\psi(1)} i}\Cal Z'_{-i, i}(w)dw \Big ) \dots 
<f_{i_r}, 
\g_{\psi(r-1)}> <f_{i_{r+1}}, \g_{\phi(1)}> \dots <f_{i_{t-1}}, 
\g_{\phi(t-r-1)}> f_{i_t}+\\
C(\Cal Z'_{-i, i}, \dots, f_{i_t})|_2(\g_1-1) 
\dots (\g_{t-2}-1).
\endmultline
\tag 3.39
$$
where $C(\Cal Z'_{-i, i}, f_{i_3}, \dots)$
is $C(f, f_{i_3}, \dots)$ with $f$ formally replaced by $\Cal Z'_{-i, i}$. 
Since the cancellations yielding (3.34') do not rely on $f_{i_1}$'s being 
a cusp form rather than second order modular form, the same 
cancellations imply that (3.39) equals
$f_{i_t}<f_{i_{t-1}}, \g_{t-2}>\dots 
\int_i^{\pi_{\ca_m} i} \Cal Z'_{-i, i}(w)dw$ 
as we wanted to show.

The second part of the claim is proved in the same way, once we observe 
that $C(f_{i_1}, f_{i_2}, \dots)$ does not contain integrals of
$\Cal Z'_{-i, i, \dots}$ because of the way the indices $j_1, j_2, \dots$
in the expression of $C(f_{i_1}, f_{i_2}, \dots)$ 
are obtained and the fact that $i_{r+1}, \dots, i_{t-1}<0$.
Therefore the inductive hypothesis can be applied.


Therefore, for all $\g_2, \dots, \g_{t-2} \in \G$,  
$$\multline
\sum 
c_{i_3, \dots, i_t} \Cal Z'_{-1, 1, \dots, i_t}|_2
(\pi_{\ca_m}-1) \dots (\g_{t-2}-1)=\\
\left ( \int_i^{\pi_{\ca_m} i} \Cal Z'_{-1, 1}(w)dw 
\right )
\sum c_{i_3, \dots, i_t} <f_{i_3}, \g_2>
\dots <f_{i_{t-1}}, \g_{t-2}> f_{i_t}
\endmultline$$ 
and since the LHS is $0$, the linear independence of the terms 
$<f_{i_3}, \g_2>
\dots <f_{i_{t-1}}, \g_{t-2}> f_{i_t}$ for $\g_2, \dots, \g_{t-2} 
\in \G$, implies the desired conclusion.
\bs

\proclaim{Theorem 3.12} Let $t \ge 1$. Then the image of 
$$\{\Cal Z_{i_1, \dots, i_t}; (i_1, \dots, i_t) \in I\}$$
under the natural projection is a basis of $S_2^t/S_2^{t-1}$.
\p
We first note that by Th. 3.1 of [CD], each $\Cal Z_{i_1, \dots, i_t}$
($(i_1, \dots, i_t) \in I$) 
is a weight $2$, $t$-th order form cusp
because it is the
product (``$0$-th Rankin-Cohen bracket") of a weight $2$, order $r<t$
cusp form and a weight $0$, order $t-r+1$ modular form 
(namely the antiderivative
of a weight $2$, order $t-r$ cusp form). 

We now show 
that our set spans $S_2^t/S_2^{t-1}$.
The claim is obvious for $t=1$. Let now $t>1$ and $F \in S_2^t.$ 
Then, for every $\g_i \in \G$,
$$F|_2(\g_1-1) \dots (\g_{t-1}-1)=\sum_i \chi_i(\g_1, \dots, \g_{t-1}) 
f_i $$
for a $\chi_i: \G^{t-1} \to \Bbb C$. 
As in the proof of Lemma 3.11, each $\chi_i$ is a group 
homomorphism but, in addition, each of them vanishes at the parabolic 
elements as a function of each $\gamma_i$. 
Repeated applications of the Eichler-Shimura 
isomorphism imply that each $\chi_i(\g_1, \dots, \g_{t-1})$ is a linear 
combination of 
$$
\prod_{j=1}^{t-1} <f_{i_j}, \g_j> \quad i_j \in \{\pm 1, \dots, \pm g\}.
$$ 
By the construction of $\Cal Z_{i_1, \dots i_t}$'s and (3.34) (with $s=t$) 
we then deduce that there is a linear combination $L$ of these functions 
satisfying,
$(F-L)|_2(\g_1-1) \dots (\g_{t-1}-1)=
\phi(\g_1, \dots, \g_{t-1})$ for some $\phi \in \Cal A$. 
We will show that $\phi \equiv 0.$

By the definition of $\Cal A$, there are $c_{i_1, \dots, \hat i_j, \hat 
i_{j+1}, \dots, i_t} \in \Bbb C$ such that $\phi(\g_1, \dots, \g_{t-1})$
equals
$$
\sum_{j=1}^{t-1}\sum c_{i_1, \dots, \hat i_j, \hat i_{j+1}, \dots, 
i_t}<f_{i_1}, \g_1> \dots 
\overline{\int_{i}^{\g_j i} f_1(w)dw} 
\int_{i}^{\g_{j+1} i} f_1(w)dw <f_{i_{j+2}}, \g_{j+2}> \dots f_{i_t}
$$
for all $\g_i \in \G$. The inner sum ranges over all $(i_1, \dots, \hat 
i_j, \hat i_{j+1}, \dots, i_t)$ such that $(i_{j+2} 
\dots, i_t) \in I$ and the hat indicates missing index. The term 
corresponding to $j=t-1$ is understood to end with 
$\overline{\int_{i}^{\g_{t-1} i} f_1(w)dw} f_1$. With (3.34), an 
induction shows that, for some $a^l_{i_1, 
\dots}, b^l_{i_{j+2}, \dots} \in \Bbb C$, the sum can be re-written as
$$\multline
\sum_{j=1}^{t-1}\sum_l \sum_{0< |i_k| \le g} a^l_{i_1, \dots, i_{j-1}} 
<f_{i_1}, 
\g_1> \dots \times \\
\overline{\int_{i}^{\g_j i} f_1(w)dw} 
\int_{i}^{\g_{j+1} i} f_1(w)dw \sum_{0< |i_k| \le g} b^l_{i_{j+2}, \dots, 
i_{t-1}} <f_{i_{j+2}}, \g_{j+2}> \dots f_{i_t}
\endmultline
\tag 3.40
$$
with 
$$\int_{i}^{\g_{j+1} i} f_1(w)dw \sum_{0< |i_k| \le g} b^l_{i_{j+2}, 
\dots, i_{t-1}} <f_{i_{j+2}}, \g_{j+2}> \dots f_{i_t}
=G|_2(\g_{j+1}-1) \dots (\g_{t-1}-1)$$
for some $G \in S_2^{t-j}.$

This implies that, if we consider $\g_1, \dots, \g_{t-2}$ fixed for 
the time being, each term in (3.40) except for that corresponding to 
$j=t-1$ is a multiple of a $G|_2(\g_{t-1}-1)$ for some $G \in S_2^2$. 
Hence, if $(F-L)|_2(\g_1-1) 
\dots (\g_{t-1}-1)=
\phi(\g_1, \dots, \g_{t-1})$, then there is a second-order form $G_1$ and 
a $\mu_{t-1} \in \Bbb C$ (which will normally depend on $\g_1, \dots, 
\g_{t-2}$) such that 
$$G_1|_2(\g_{t-1}-1)=\mu_{t-1} \overline{\int_{i \ci}^{\g_{t-1} i\ci} 
f_1(w)dw}f_1 $$
Prop. 5.2 of [DO] implies that $\mu_{t-1}=0.$ Therefore, $\phi(\g_1, 
\dots, \g_{t-1})$ equals the expression in (3.40) but with one term less.

We can continue the `descent' this way by noting each term on the 
right-hand side except for the last one is a multiple of a
$G|_2(\g_{t-2}-1)(\g_{t-1}-1)$ for some $G \in S_2^3$. 
Therefore, there is a $G_2 \in S_2^3$ such that
$$G_2|_2(\g_{t-2}-1)(\g_{t-1}-1)=\sum c_{i_t}
\overline{\int_{i \ci}^{\g_{t-2} i\ci} f_1(w)dw} \int_{i 
\ci}^{\g_{t-1} i\ci} f_1(w)dw f_{i_t} $$
for some $c_{i_t} \in \Bbb C$
and this, from Lemma 3.11, implies that all $c_{i_t}$ vanish. Continuing 
this way we deduce that $\phi \equiv 0$.

We will finally show that the $\Cal Z$'s are linearly independent modulo 
$S_2^{t-1}$. 
Suppose that 
$$\sum_{(i_1, \dots, i_t) \in I}\lambda_{i_1, \dots, i_t} \Cal Z_{i_1, 
\dots, i_t}=0. $$ 
Then
$$\sum_{(i_1, \dots, i_t) \in I}\lambda_{i_1, \dots, i_t} \Cal Z_{i_1, 
\dots, i_t}|_2(\g_1-1)\dots(\g_{t-1}-1)=0$$
and by $(3.34)$ (with $s=t$), we obtain
$$\sum_{(i_1, \dots, i_t) \in I}\lambda_{i_1, \dots, i_t} 
<f_{i_1}, \g_1>\dots <f_{i_{t-1}}, \g_{t-1}> f_{i_t} +\phi(\g_1, \dots, 
\g_{t-1})=0,$$
for some $\phi \in \Cal A$. Since the set of all $\chi_i$'s is 
linearly independent and $\phi$ is a linear combination of 
$<f_{i_1}, \g_1>\dots <f_{i_{t-1}}, \g_{t-1}> f_{i_t}$'s
with $(i_1, \dots, i_t) \not \in I$, we deduce that $\lambda_{i_1, \dots, 
i_t}=0$, for all $(i_1, \dots, i_t) \in I$, i.e. the $\Cal Z$'s are 
linearly independent. \bs

\proclaim{Corollary 3.13} Let $g$ be the genus of $\G$. The 
dimension of $S_2^t(\G)/S_2^{t-1}(\G)$ ($t \ge 1$) is $0$ if $g=0$, and
$$\frac{1}{2} \Big ( (g+\sqrt{g^2-1})^t+(g-\sqrt{g^2-1})^t\Big )$$
otherwise.
\p 
According to the theorem, to prove the corollary in the case $g \ne 0$, it 
suffices to compute the cardinality $a_t$ of $I$. If $b_t$ is the 
cardinality of $I'$, then, $b_t-a_t$ 
is the number of vectors in $I'$ with $i_{t}<0$. In particular, $i_t \ne 
1$ and hence, $b_t-a_t=gb_{t-1}$. On the other hand,
$b_t=2gb_{t-1}-b_{t-2}$. ($gb_{t-1}$ elements of $I'$ have $i_t<0$ and 
$gb_{t-1}-b_{t-2}$ have $i_t>0$ with $(i_{t-1}, i_t) \ne (-1, 1)$).
The solution of this recursive relation with $b_1=2g, b_2=4g^2-1$ is
$$b_t=\frac{1}{2 \sqrt{g^2-1}}
\Big ((g+\sqrt{g^2-1})^{t+1}-(g-\sqrt{g^2-1})^{t+1}\Big )$$
when $g \ne 1$ and $b_t=t+1$, when $g=1$. 
The conclusion then follows from $b_t-a_t=gb_{t-1}.$

The proof of the corollary in the case $g=0$ follows by the observation 
that, for $F \in S_2^t(\G)$,
$$F|_2(\g_1-1) \dots (\g_{t-1}-1) \in S_2(\G). \tag 3.41$$
Since $S_2(\G)=\{0\}$ when $g=0$, this implies that $F \in S_2^{t-1}(\G)$
and hence $S_2^t(\G)/S_2^{t-1}(\G)=\{0\}$.
\bs
\bf Examples \rm (with $g \ne 0$):  $$\matrix t & \, \, & \dim 
(S_2^t(\G)/S_2^{t-1}(\G)) \\  
\, \\
1 & \, \, & g  \\
2 & \, \, & 2g^2-1 \\
3 & \, \, & 4g^3-3g \\
4 & \, \, & 8g^4-8g^2+1 \\
5 & \, \, &16g^5-20g^3+5g \\
\endmatrix$$

$$\text{\bf 4. Higher Weights}$$

Our construction of the base of $S_k^t$ for $k>2$ and $t>1$ relies on the 
base for $S_2^{t-1}$ we have just defined and it parallels the process 
we employed to construct the latter. We shall maintain the notation 
developed in the previous sections.

First of all we note that Theorem 3.4 can be used just as well to 
construct forms of higher weight. Specifically, for every $k>2$ and $f 
\in A_t$ we set
$$\align Y_{\ca m}(z, s; \bar f)& :=y^{-k/2}Z_{\ca m}(z, s+k/2, 1, k; \bar 
f)\\
&=\sum_{\g \in \G_\ca \backslash \G} \overline{\Big (\int_{\ca}^{\g \ca}
f(w) dw \Big )}
\Im(\sa^{-1} \g z)^s e(m \sa^{-1} \g z) j(\sa^{-1} \g, z)^{-k}.
\endalign
$$
According to Theorem 3.4, it is analytic for Re$(s)>1-k/2-\delta_{\G}$
and it is easy to see that
$$\multline
Y_{\ca m}(\cdot, 0, \bar F_{i_1, \dots, i_{t-1}})|_k(\g-1)=\\
\Big (\overline{\int_{\ca}^{\g^{-1} \ca} 
F_{i_1, \dots, i_{t-1}}}\Big ) P_{\ca m}(\cdot)_k+
\sum_{r=1}^{t-2} \Big (\overline{\int_{\ca}^{\g^{-1} \ca} 
F_{i_1, \dots, i_r}}\Big )Y_{\ca m}(\cdot, 0, 
\bar F_{i_{r+1}, \dots, i_{t-1}})
\endmultline
\tag 4.1
$$
Here we set
$$P_{\ca m}(z)_k:=y^{-1}U_{\ca m}(z, 1, k),$$
the Poincar\'e series of weight $k$. There is a set of positive 
integers $M$ such that $\{P_{\ca m}(\cdot)_k\}_{m \in M}$ is a basis of 
$S_k$.

Now, for Re$(s)$ large we have:
$$
\frac{d}{d\bar z}Y_{\ca m}(z, s; \bar F_{i_1, \dots, i_{t-1}})=
\frac{is}{2y^{1+k/2}}Z_{\ca m}(z, s+k/2, 1, k-2; \bar F_{i_1, \dots, 
i_{t-1}})
$$
By Theorem 3.4., $Z_{\ca m}(z, s+k/2, 1, k-2; \bar F_{i_1, \dots,
i_{t-1}})$ is holomorphic at $s=0$, when $k>2$ and hence
$Y_{\ca m}(z, 0; \bar F_{i_1, \dots, i_{t-1}})$ is holomorphic in $z$.
Thus we set $\Cal Y_m:= P_{\ca m}(\cdot)_k$ and
$$\Cal Y_{-i_1, \dots, -i_{t-1}; m}=(-1)^{t-1}Y_{\ca m}(z, 0; \bar F_{i_1, 
\dots, i_{t-1}}).$$
With (4.1) we note that
$$\Cal Y_{-i_1, \dots, -i_{t-1}; m}|_k(\g_1-1) \dots (\g_{t-1}-1)=
\overline{\int_{\ca}^{\g_{1} \ca} f_{i_1}} \dots
\overline{\int_{\ca}^{\g_{t-1} \ca} f_{i_{t-1}}}
P_{\ca m}(\cdot)_k.$$
In exactly the same way as Th. 3.9, we can show that $\Cal Y_{i_1, \dots, 
i_t; m} \in S_k^{t}$ for $i_j \in \{-1, \dots, -g\}$.
  
We can now construct a basis for $S_k^{t+1}$ by extending the definition 
of $
\Cal Y_{i_1, \dots, i_t; m}$ to all $(i_1, \dots, i_t) \in I'$ and $m>0$:
First, with Lemma 3.10, for $(i_1, \dots, i_t) \in I$
$$\multline
[\Cal Y_{m} (z)\int_{i}^z \Cal Z_{i_1, \dots, i_t}(w)dw]|_k(\g_1-1) 
\dots (\g_{t}-1)=\\
P_{\ca m}(\cdot)_k <f_{i_1}, \g_1> \dots <f_{i_{t-1}}, \g_{t-1}> 
\int_{i}^{\g_t i} f_{i_t}(w)dw +\phi(\g_1, \dots, \g_{t})
\endmultline
$$
for some $\phi \in \Cal A_k.$
Hence, for 
$(i_1, \dots, i_{t}) \in I$ 
we set 
$$\Cal Y_{i_1, \dots, i_t; m}:=\Cal Y_{m}(z) 
\int_i^z \Cal Z_{i_1, \dots, i_{t}}(w)dw.$$

Next, if $i_1<0$, and $(i_2, \dots, i_t) \in I$,
$$\multline
[\Cal Y_{i_1; m}(z) \int_{i}^z \Cal Z_{i_2, \dots, 
i_t}(w)dw]|_k(\g_1-1) \dots (\g_{t}-1)=\\
\sum_{(\phi, \psi) \in \Cal{S}_{2, t+1}} $$
$$\overline{\int_{i}^{\g_{\phi(1)} i}f_{i_1}(w)dw} 
P_{\ca m}(\cdot)_k \prod_{j=2}^{t} <f_{i_j}, \gamma_{\psi(j)}> 
\int_{i}^{\g_{\psi(t)} i} f_{i_{t}}(w)dw+ \phi(\g_1, \dots, 
\g_{t})
\endmultline
$$
for some $\phi \in \Cal A_k$.
Therefore, as in Section 3.4, there is a linear combination of 
$\Cal Y_{m} (z)\int_{i}^z \Cal Z_{j_1, \dots, j_t}(w)dw$'s 
denoted by $A_{i_1, \dots, i_t; m}$, such that 
$$\multline
[\Cal Y_{i_1, m} (z)\int_{i}^z \Cal Z_{i_2, \dots, 
i_t}(w)dw-A_{i_1, \dots, i_t; m}]|_k(\g_1-1) \dots (\g_t-1)=\\
\overline{\int_{i}^{\g_t i} f_{i_1}(w)dw} 
P_{\ca m}(\cdot)_k <f_{i_2}, \gamma_1> \dots <f_{i_{t-1}}, \g_{t-2}> 
\int_{i}^{\g_{t-1} i} f_{i_t}(w)dw+ \phi(\g_1, \dots, \g_{t})
\endmultline
$$
for some $\phi \in \Cal A_k$.
Hence for $(i_{1}, \dots, i_{t-1}) \in I$ and $i_{t}<0$ 
we can set 
$$\Cal Y_{i_1, \dots, i_t; m}(z):=\Cal Y_{i_{t}; m}(z) \int_{i}^z \Cal 
Z_{i_1, \dots, i_{t-1}}(w)dw-A_{i_{t}, i_1, \dots, i_{t-1}; m}.$$

Further, if $i_1, i_2 <0$ and $(i_3, \dots, i_t) \in I$,
$$\multline
[\Cal Y_{i_1, i_2; m} (z)\int_i^z \Cal Z_{i_3, \dots, 
i_t}(w)dw]|_k(\g_1-1) \dots (\g_{t}-1)=\\
\sum_{(\phi, \psi) \in \Cal{S}_{3, t+1}} 
\overline{\int_{i}^{\g_{\phi(1)} i} f_{i_1}(w)dw}
\overline{\int_{i}^{\g_{\phi(2)} i} f_{i_2}(w)dw}
P_{\ca m}(\cdot)_k\prod_{j=3}^{t} <f_{i_j}, \gamma_{\psi(j)}> 
\int_{i}^{\g_{\psi(t)} i} f_{i_{t}}(w)dw+ \\
\phi(\g_1, \dots, \g_{t})\endmultline
$$
for some $\phi \in \Cal A_k$.
Therefore, there is a linear combination of 
terms $\Cal Y_{m}(z) \int_{i}^z \Cal Z_{j_1, \dots, 
j_t}(w)dw$'s and 
$\Cal Y_{j_1; m}(z) \int_{i}^z \Cal Z_{j_2, \dots, j_t}(w)dw$'s 
denoted by $B_{i_1, \dots, i_t; m}$, such that 
$$[\Cal Y_{i_1, i_2; m}(z) \int_{i}^z \Cal Z_{i_3, \dots, 
i_t}(w)dw-B_{i_1, \dots, i_t; m}]|_k(\g_1-1) \dots (\g_{t}-1)=$$
$$\overline{\int_{i}^{\g_{t-1} i} f_{i_1}(w)dw}
\overline{\int_{i}^{\g_t i} f_{i_2}(w)dw}
P_{\ca m}(\cdot)_k<f_{i_3}, \gamma_1> \dots \int_{i}^{\g_{t-2} i} 
f_{i_t}(w)dw+ \phi(\g_1, \dots, \g_{t})$$
for some $\phi \in \Cal A_k$.
Hence for $(i_1, \dots, i_{t-2}) \in I$ and 
$i_{t-2}, i_{t-1}<0$ we can set
$$\Cal Y_{i_1, \dots, i_t; m}(z):=\Cal Y_{i_{t-1}, i_{t}; m}(z) \int_i^z 
\Cal Z_{i_1, \dots, i_{t-2}}(w)dw-B_{i_{t-1}, i_{t}, i_1, 
\dots, i_{t-2}}.$$

Continuing this way we cover the entire $I'$ and hence we construct, for 
all $(i_1, \dots, i_t) \in I'$, 
functions $\Cal Y_{i_1, \dots i_t; m}$ such that 
$$\Cal Y_{i_1, \dots, i_t; m}|_k(\g_1-1) \dots (\g_t-1)=
<f_{i_1}, \g_1>\dots <f_{i_t}, \g_t> P_{\ca m}(\cdot)_k 
+\phi(\g_1, \dots, \g_t) 
\tag 4.2$$
for some $\phi \in \Cal A_k$.

\proclaim{Theorem 4.1}
If $k\geqslant 4$ and $g$ is the genus of $\G$, then the image of the 
set $$\{\Cal Y_{i_1, \dots, i_t; m}; (i_1, \dots, i_{t}) \in I'; m \in M 
\}$$
under the natural projection is a basis of $S_k^{t+1}(\G)/S_k^{t}(\G)$. 
Therefore the dimension of $S_k^{t+1}(\G)/S_k^{t}(\G)$ is $0$ if $g=0$, 
$(t+1)\dim(S_k(\G))$ if $g=1$, and
$$\frac{\dim(S_k(\G))}{2\sqrt{g^2-1}} 
\Big ( 
(g+\sqrt{g^2-1})^{t+1}-(g-\sqrt{g^2-1})^{t+1}\Big )$$
otherwise.
\p
Each $\Cal Y_{i_1, \dots, i_t; m}$ is a $(t+1)$-th order cusp form 
as the product of an order $r<t+1$ cusp form and an order $t+2-r$ modular 
form.

Now, since the argument does not depend on the weight, we can show exactly 
as in the proof of Th. 3.12, that if $F \in S_k^{t+1}$, then
$$F|_k(\g_1-1) \dots (\g_t-1)=
\sum_{i_1, \dots, i_{t+1}} \lambda_{i_1, \dots, i_{t}} 
<f_{i_1}, \g_1> \dots <f_{i_{t}}, \g_{t}> F_{i_{t+1}}. \tag 4.3$$
Here $F_i$ ranges over a basis of $S_k$.
By (4.2) we then conclude that there is a linear combination $L$ of 
functions in the set under consideration such that 
$$(F-L)|_k(\g_1-1)\dots (\g_t-1)=
\phi(\g_1, \dots, \g_t)$$
for some $\phi \in \Cal A_k$. It is easy to see that Lemma 3.11 holds in 
weights $>2$. The only adjustment required in the proof is that 
the construction of the analogue of $\Cal Z'$ is based on 
$$\Cal Y_{i_1, \dots, i_r; m} \int_i^z \Cal Z'_{i_{r+1}, \dots, i_t}$$
instead of  
$\Cal Z'_{i_1, \dots, i_r} \int_i^z \Cal Z'_{i_{r+1}, \dots, i_t}.$
We deduce that $\phi \equiv 0$ and hence, our set spans 
$S_k^{t+1}/S_k^{t}$.

The proof that the set is linearly independent is deduced directly from 
(4.2), exactly as in 
Th. 3.12.

The formula for the dimension in the case $g >0$ is deduced by the formula 
for $b_t:=\# I'$ established in Cor. 3.13. In the case $g=0$, the 
dimension is $0$ because of (4.3).
\bs

\vskip -2mm

\flushpar
\bf Examples \rm (with $g \ne 0$): 

$$\matrix t & \, \, & \dim
(S_k^t(\G)/S_k^{t-1}(\G)) \\
\, \\
2 & \, \, & 2g\dim(S_k(\G))  \\
3 & \, \, & (4g^2-1)\dim(S_k(\G)) \\
4 & \, \, & (8g^3-4g)\dim(S_k(\G)) \\
5 & \, \, & (16g^4-12g^2+1)\dim(S_k(\G)) \\
\endmatrix$$

\bigskip

\it Acknowledgments. \rm The authors thank Paul Gunnells and Cormac 
O'Sullivan for several helpful comments. We are also indebted to the 
referee for a careful reading of the manuscript and for many useful 
suggestions.

\end